\documentclass[]{siamart220329}

\usepackage{lipsum}
\usepackage{amsfonts}
\usepackage{graphicx}
\usepackage{epstopdf}
\usepackage{algorithmic}
\ifpdf
  \DeclareGraphicsExtensions{.eps,.pdf,.jpg,.jpg}
\else
  \DeclareGraphicsExtensions{.eps}
\fi

\newsiamremark{remark}{Remark}
\newsiamremark{hypothesis}{Hypothesis}
\crefname{hypothesis}{Hypothesis}{Hypotheses}
\newsiamthm{claim}{Claim}

\title{Realistic pattern formations on surfaces {by adding  arbitrary roughness}\thanks{Submitted revised version to the editors: {October, 2023.}
\funding{This work was funded by the Hong Kong Research Grant Council GRF Grants (12303818,12301419,12301520,12301021), {a National Youth Science Foundation of China  (12201449)}, and the financial support of NSERC Canada (RGPIN-2022-03302).}}}

\author{Siqing LI\thanks{{College of Mathematics}, Taiyuan University of Technology, Shanxi, China
  (\email{lisiqing@tyut.edu.cn}).}
\and Leevan LING\thanks{Department of  Mathematics, Hong Kong Baptist University, Kowloon Tong, Hong Kong
  (\email{lling@hkbu.edu.hk}).}
\and Steven J. RUUTH\thanks{Department of  Mathematics, Simon Fraser University, Burnaby, British Columbia, Canada V5A1S6
  (\email{sruuth@sfu.ca}).}
\and Xuemeng Wang\thanks{Department of Mathematics, University of British Columbia,   Vancouver, Canada, and Department of  Mathematics, Hong Kong Baptist University, Kowloon Tong, Hong Kong
  (\email{17251109@life.hkbu.edu.hk}).}
}

\usepackage{amsopn}

\usepackage[percent]{overpic}
\usepackage{graphicx,graphics,epsfig}
\usepackage{amsmath,amssymb,dsfont}
\usepackage{mathrsfs}
\usepackage{arydshln}
\usepackage{url,cite}
\usepackage{enumerate}
\usepackage{subfigure}
\usepackage{algorithm}
\usepackage{amsmath}
\usepackage{graphicx}
\usepackage{booktabs}
\usepackage{threeparttable}
\usepackage{bm}
\usepackage{multirow}
\usepackage{multicol}

\usepackage[normalem]{ulem}

\def\eref#1{{\rm (\ref{#1})}}

\def\qed{~\relax\ifmmode\hskip2em \Box
 \else\unskip\nobreak\hskip1em \hfill$\Box$
 \fi \newline}

\def\D{\mathbb{D}}

\def\R{\mathbb{R}}

\def\cala{\mathcal{A}}

\def\calf{\mathcal{F}}

\def\cali{\mathcal{I}}

\def\calm{\M}

\def\cals{\mathcal{S}}

\def\calu{\mathcal{U}}

\def\bigdot{\boldsymbol{\cdot}}

\def\D{{\mathcal{D}}}
\def\M{{\mathcal{M}}}
\def\MAmp{\delta_{\M}}
\def\bx{\boldsymbol{\xi}}

\usepackage{newunicodechar}

\begin{document}

\maketitle

\begin{abstract}
We are interested in generating surfaces with arbitrary roughness and forming patterns on the surfaces. Two methods are applied to construct rough surfaces. In the first method, some superposition of wave functions with random frequencies and angles of propagation are used to get periodic rough surfaces with analytic parametric equations. The amplitude of such surfaces is also an important variable in the provided eigenvalue analysis for the Laplace-Beltrami operator and in the generation of pattern formation. Numerical experiments show that the patterns become irregular as the amplitude and frequency of the rough surface increase. For the sake of easy generalization to closed manifolds, {we propose a second construction method for rough surfaces,  which uses random nodal values and discretized heat filters.} We provide numerical evidence that both surface {construction methods} yield comparable patterns to those {observed} in real-life animals.
\end{abstract}

\begin{keywords}
Laplace-Beltrami operator, reaction-diffusion system, random surfaces, Turing pattern
\end{keywords}

\begin{MSCcodes}
65M06,  
35K57 
\end{MSCcodes}

\section{Introduction}

\begin{figure}
	\centering
	\setlength{\abovecaptionskip}{0pt}
	\setlength{\belowcaptionskip}{0pt}
	\begin{tabular}{cc}
	\begin{overpic}[width=0.45\textwidth]{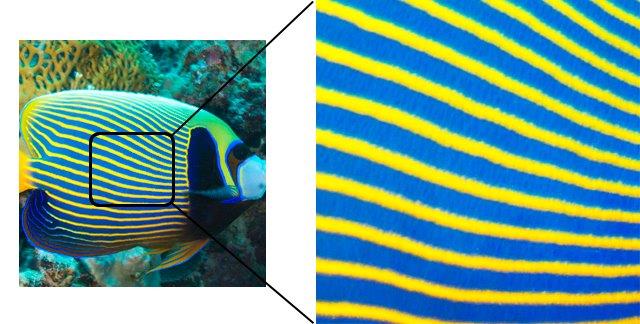}
		\put(0,0) {\scriptsize (a) Emperor angelfish}
	\end{overpic}
	&
	\begin{overpic}[width=0.45\textwidth]{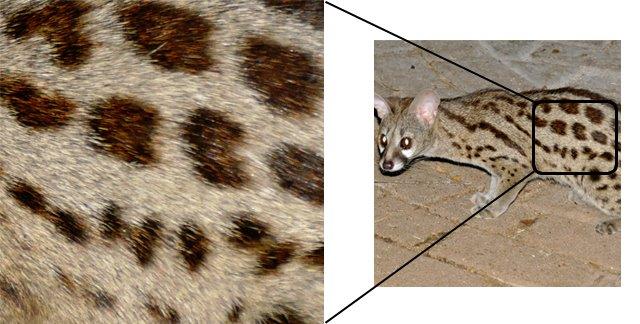}
		\put(65,0) {\scriptsize (b) Genet}
	\end{overpic}
\\
	\begin{overpic}[width=0.45\textwidth]{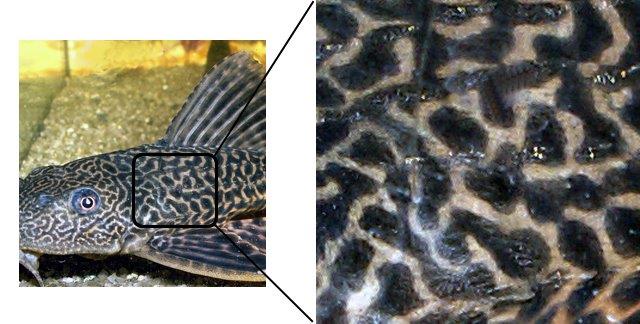}
		\put(0,0) {\scriptsize (c) Plecostomi  }
	\end{overpic}
	&
	\begin{overpic}[width=0.45\textwidth]{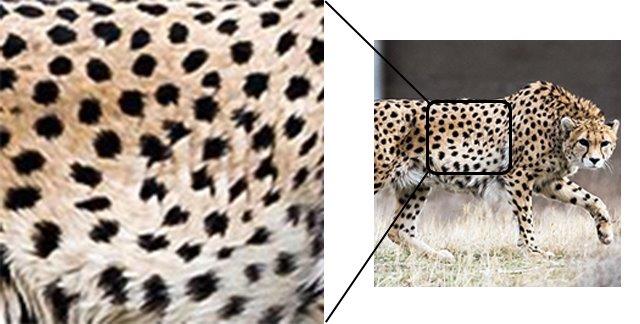}
			\put(65,0) {\scriptsize (d) Cheetah}
	\end{overpic}
	\end{tabular}
  \caption{{Real-life examples of skin patterns.}
  }\label{fig:motivation}
\end{figure}

Pattern formation by reaction-diffusion systems has been an intensively studied field for decades.
In 1952, Turing proposed the idea of diffusion-driven instability  \cite{chemicalbasis}, wherein simple mechanisms evolve from a homogeneous state into spatial heterogeneous patterns.
In recent years, various application areas have been actively developed, including vegetation {patterns}, plant root hair initiation, flock formation and boundary drop configurations \cite{doi:10.1137/21M1390141, doi:10.1137/17M1120932, doi:10.1137/17M113856X, doi:10.1137/18M1196996}. Mechanisms of pattern generation under different types of transport
or domain size have
been {explored} \cite{doi:10.1137/16M1061205, doi:10.1137/17M1138571}. In theoretical biology, reaction-diffusion systems provide a relatively generic and concise approach for generating animal skin patterns \cite{coat,mathematicalaspects}, as well as regenerative processes of organisms \cite{GM}. In particular, reaction-diffusion systems are a well-accepted class of models for multiple pigmentation processes, including emperor  angelfish \cite{Albert2014File},
genets \cite{Category},
plecostomi \cite{Contributors2021Plecostomus},
cheetah  \cite{2015undefined} patterns, which are shown in Figure~\ref{fig:motivation}, as well as
marine angelfish \cite{Kondo1995ARW},
zebra fish\cite{asai_zebrash_1999}, and various mammal skin patterns \cite{Murray1980APF}.
Although animal skin pattern through reaction-diffusion systems has been explored in many researches, see the second row of Figure~\ref{fig:ex31} for typical examples,
very few {studies} have {considered} the skin texture. In reality, the skin surfaces are often rugged, {resulting in patterns that are not strictly characterized by regular spots or stripes. }
{Therefore, combining surface properties with reaction-diffusion systems can provide an opportunity for improving generation of real patterns. }
In \cite{AAMM-12-1327,Chu_2021}, parameter functions instead of constants were used in PDEs to generate nonuniform and complex  patterns {that closely resemble real-life patterns.}
%

Characterizing random rough surfaces can be {achieved} by a variety  methods.  These include the reference parameter method, motif method, fractal method, watershed method, and wavelet method, etc \cite{shi_effect_2019}. For three-dimensional surface topography, there are a number of parameters involved \cite{gadelmawla_roughness_2002}.
To generate random rough surfaces, {two commonly used methods are the Fast Fourier Transform (FFT)\cite{pawlus2020review} and digital filter \cite{hu1992simulation} methods. These methods  approximate the auto-correlation function for surfaces with wavelengths. }
In \cite{huang2021isogeometric}, the authors proposed an approach that applies the covariance function and Karhunen-Lo\`{e}ve expansion to generate random surfaces. {Alternatively, rough surfaces can also be generated from spatial frequencies via the FFT method }by summing up trigonometric functions \cite{sjodin_how_nodate}.


Moving from a flat geometry to rough surfaces generally leads to surface dependent differential operators. Surface PDEs have been extensively studied in recent years, and {various techniques exist }for their analysis and computation.
The approaches  for surface PDEs can be {categorized as}
extrinsic, embedding, and intrinsic methods.
Extrinsic methods solve PDEs on the manifold directly.
Some extrinsic methods depend on surface mesh construction.
Examples of such methods include finite difference \cite{ruuth_simple_2008}, finite element \cite{dziuk_finite_2013,olshanskii_finite_2009,bertalmio_variational_2001}, and finite volume \cite{du_finite_2005,du_voronoi-based_2003} methods. {Additionally, there are }numerical methods that do not require meshes, {known as} meshfree methods. Meshfree methods for surface PDEs, such as radial basis function (RBF) methods \cite{rbf, Chen+Ling-Extrmeshcollmeth:20, Chen-KernMeshCollMeth:19} and the meshfree generalized finite difference method \cite{Tang+FuETAL-locaextrcollmeth:21
,SUCHDE20192789},  have the advantage of avoiding mesh construction.
On the other hand, embedding methods \cite{Cheung+Ling-Kernembemethconv:18,Maerz+Macdonald-CALCSURFWITHGENE:12} formulate and solve PDEs on a band around the manifold.
{In this paper}, intrinsic methods \cite{Charette+MacdonaldETAL-Pattformslowflat:20,Sun-ConvIterQuasPeri:22} are used, {which work on parameter spaces of (local) parametrization.}

Our paper focuses on the construction of surfaces with arbitrary roughness and the generation of patterns by solving PDEs on surfaces.  In Section \ref{sec:introductionroughsurface}, the rough surface $\calm$ with parametric equations will be formulated in detail. We consider periodic rough surfaces $\calm$ characterized with spatial frequency, which {are accompanied by} analytic parametric equations. The eigenvalue analysis for the Laplace-Beltrami operator is  provided. Along with the construction of rough surfaces, {Section~\ref{sec:FDM} reviews the}
intrinsic concept of the surface heat equation, which is a PDE involving the Laplace-Beltrami operator defined by  Riemannian metric. We apply the finite difference method to solve the heat equations, {which} is chosen because of its simplicity and effectiveness on the class of problem domains considered.  In Section~\ref{sec:num}, we {employ} the finite difference method to solve pattern formation PDEs on rough surfaces $\calm$. Our aim is to study the effect of roughness on pattern formation using the simplest setup. Numerical experiments for the approximation of pattern formation on rough surfaces $\calm$ with parametric equations are demonstrated. Results for PDEs with constant parameters on rough surfaces with various roughness are similar to the real-life patterns displayed in Figure~\ref{fig:motivation}. The rough surface pattern formation models here are not yet extendable to add roughness to other manifolds; e.g., we cannot yet create red-blood cells with rough surfaces.
In Section~\ref{sec:RRS}, we {present an alternative} construction method of rough surfaces $\cals$  by random discrete data and discretized heat filters.  After updating our numerical schemes for $\calm$ to work on $\cals$,  we {observe} that surface types $\calm$ and $\cals$ yield comparable patterns and animal coats. Furthermore, the proposed surface construction technique can be {readily} extended to {introduce}  ``roughness'' to more general manifolds. {For instance, it can be applied} to generate red-blood cells with rough surfaces.



\section{Rough surfaces $\calm$ with analytic parametric equations} \label{sec:introductionroughsurface}

{We} consider $\mathcal{C}^k$--smooth ($k \geq 2$), codimension 1, and  periodic  Riemannian  surfaces{, which can be described as}
\begin{equation}\label{eq:M}
    \M=\big\{(x,y,z)\in\mathbb{R}^3 \,:\, z = z(x,y) \;\;\mbox{for } (x,y)\in V \subset \R^2\big\} ,
\end{equation}
for some $\mathcal{C}^k$ function $z(\cdot,\cdot)$ defined on a global parameter space $V\subset\R^2$.
A corresponding parametric representation of $\M$ is given by
\[
    \vec r(x,y) = \big[x,y,z(x,y)\big]^T \in \M, \quad (x,y) \in V.
\]
For convenience, {we will interchangeably use $(x,y)$ and $(x_1,x_2)$ to refer to the parameters.}
Recall that the first fundamental form $G:V\to \R^{2\times2}$ of $\calm$ is defined by
\[
G(x,y) = \begin{bmatrix}
     g_{11} & g_{12} \\
      g_{21} & g_{22}
\end{bmatrix}(x,y)
\quad \mbox{ where }\quad
g_{ij} (x,y) = \frac{\partial \vec  r(x,y) }{\partial x_i} \boldsymbol{\cdot} \frac{\partial \vec  r(x,y) }{\partial x_j}{;}
\]
see \cite{dziuk_finite_2013} for a review. In the case of a surface function $\calm$ as in \eref{eq:M}, we have
\begin{equation}
\label{eq:G}
G(x,y) = \begin{bmatrix}
     g_{11} & g_{12} \\
      g_{21} & g_{22}
\end{bmatrix}(x,y)
=\begin{bmatrix}
        1+z_x^2 & z_x z_y \\
        z_x z_y & 1+z_y^2 \\
\end{bmatrix},
\end{equation}
whose determinant and inverse are
\begin{equation}
\label{eq:g}
    g(x,y) :=\det(G)(x,y)  =1+z_x^2+z_y^2,
\end{equation}
and
\begin{equation}
\label{eq:iG}
    G^{-1}(x,y) =
\begin{bmatrix}
    g^{11}  & g^{12} \\
   g^{21} &    g^{22}
\end{bmatrix}(x,y)  =\frac{1}{1+z_x^2+z_y^2} \left( \begin{array}{c c}
        1+z_y^2 & -z_x z_y \\
        -z_x z_y & 1+z_x^2 \\
      \end{array}
     \right) .
\end{equation}
Using \eref{eq:g}--\eref{eq:iG}, the Laplace-Beltrami operator on $\M$ is given by
\begin{equation}\label{eq:LapBeltrami}
 \Delta_{\M} f= \frac{1}{\sqrt{g}} \sum_{1\leq i,j \leq 2} \frac{\partial}{\partial_i}
 \left( \sqrt{g} g^{ij} \frac{\partial}{\partial_j }f\right) =:
  {\frac1{\sqrt{g}}}\nabla \bigdot \big( {\cala} \nabla f \big),
\end{equation}
for any $\mathcal{C}^2$ function $f:\calm\to\R$. 

%

{Now, our focus shifts to random rough surfaces $\M \subseteq \mathbb{R}^3$ described in the form of \eref{eq:M}, where the surface function {$z:\cali^2=[-L,L]^2\to\R$} is a superposition of elementary waves \cite{sjodin_how_nodate,fractalimage}.} The function $z$ is stochastically determined by
\begin{equation}\label{eq:roughsurface1}
\begin{aligned}
z(x,y)=\sum_{m=-M}^{M} \sum_{n=-N}^{N}  a_{m,n}\cos\Big(2\pi(mx+ny)+\phi_{m,n}\Big),
\end{aligned}
\end{equation}
for some random variables $a_{m,n}$ and $\phi_{m,n}$.
Similar to a Fourier series expansion, $z(x,y)$ in \eref{eq:roughsurface1}  is constructed by trigonometric functions, {where}
$m,\ n$ correspond to spatial frequencies on the $x$ and $y$ axes{,} respectively. As in \cite{sjodin_how_nodate}, the spatial frequencies $m$ and $n$ allow values taken up to maximum integers $M$ and $N${, respectively. This} corresponds to a high frequency cut off.

\begin{figure}
	\centering
	\setlength{\abovecaptionskip}{0pt}
	\setlength{\belowcaptionskip}{0pt}
	\begin{tabular}{cc}
	\begin{overpic}[width=0.42\textwidth]{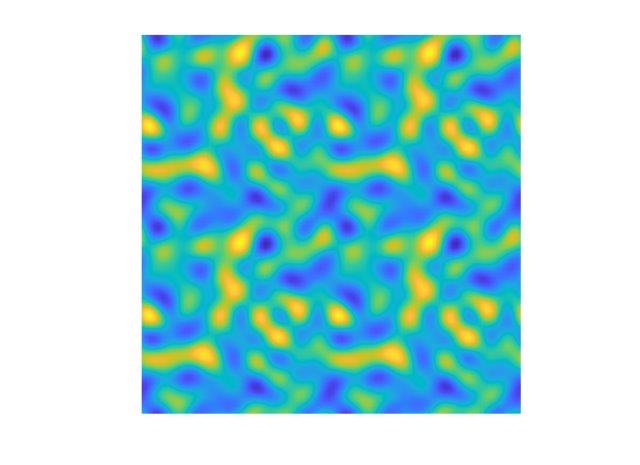}
		\put(40,70) { {  {$\M$}}}
		\put(35,3) {\scriptsize (a) $M=5=N$}
	    \put(8,15){\rotatebox{91}{ {Lower roughness} }}
	\end{overpic}
	&
	\begin{overpic}[width=0.42\textwidth]{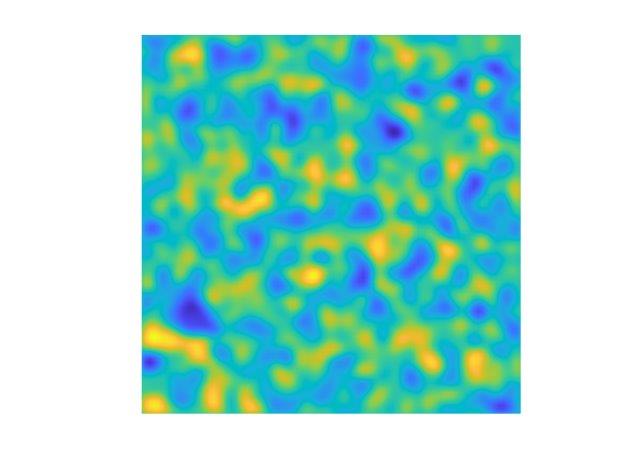}
        \put(40,70) { {  {$\cals$}}}
		\put(35,3) {\scriptsize (d) $M=5=N$ }
	\end{overpic}
	\\
	\begin{overpic}[width=0.42\textwidth]{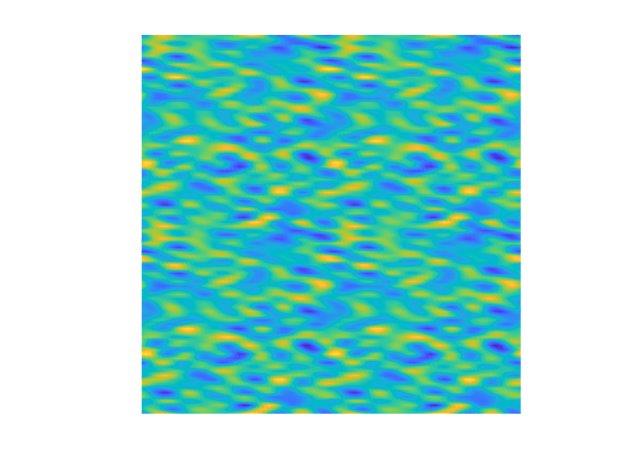}
		\put(30,3) {\scriptsize (b) $M=5,\,N=15$}
        \put(8,10){\rotatebox{91}{ $\xleftarrow{\quad\text{Increasing $M,N$}\quad}$ }}
	\end{overpic}
	&
	\begin{overpic}[width=0.42\textwidth]{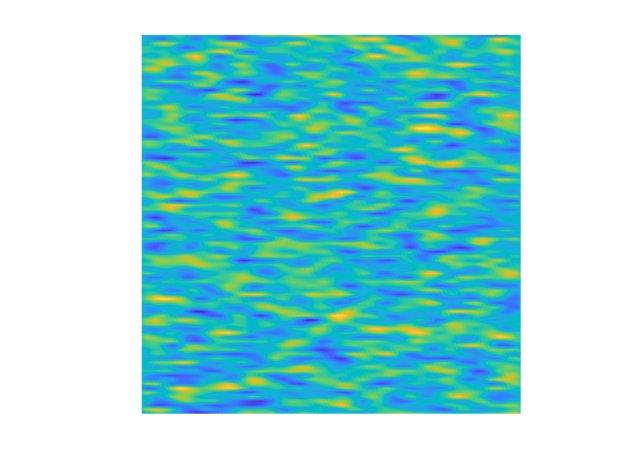}
		\put(30,3) {\scriptsize (e) $M=5,\,N=15$}
	\end{overpic}
	\\
	\begin{overpic}[width=0.42\textwidth]{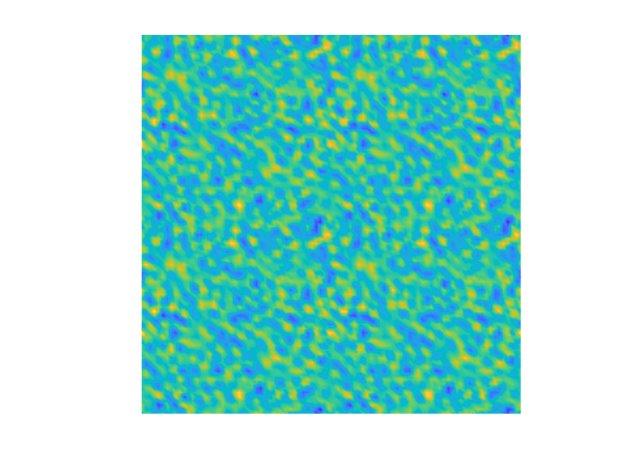}
		\put(35,3) {\scriptsize (c) $M=15=N$}
	    \put(8,15){\rotatebox{91}{ {Higher roughness} }}
	\end{overpic}
    &
	\begin{overpic}[width=0.42\textwidth]{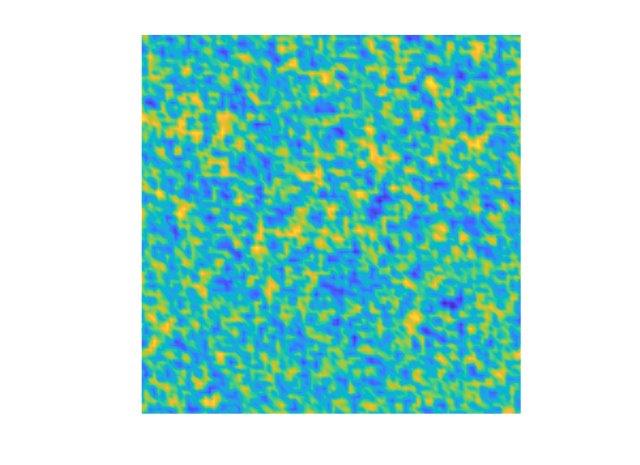}
		\put(35,3) {\scriptsize (f) $M=15=N$}
	\end{overpic}
	\end{tabular}
	\caption{{A bird's-eye view depicting} random rough surfaces $\M$ {(left) by the construction method } from Section \ref{sec:introductionroughsurface}  and $\mathcal{S}$ {(right) by the construction method} from Section \ref{sec:RRS} .
}
	\label{fig:roughsurface}
	\end{figure}
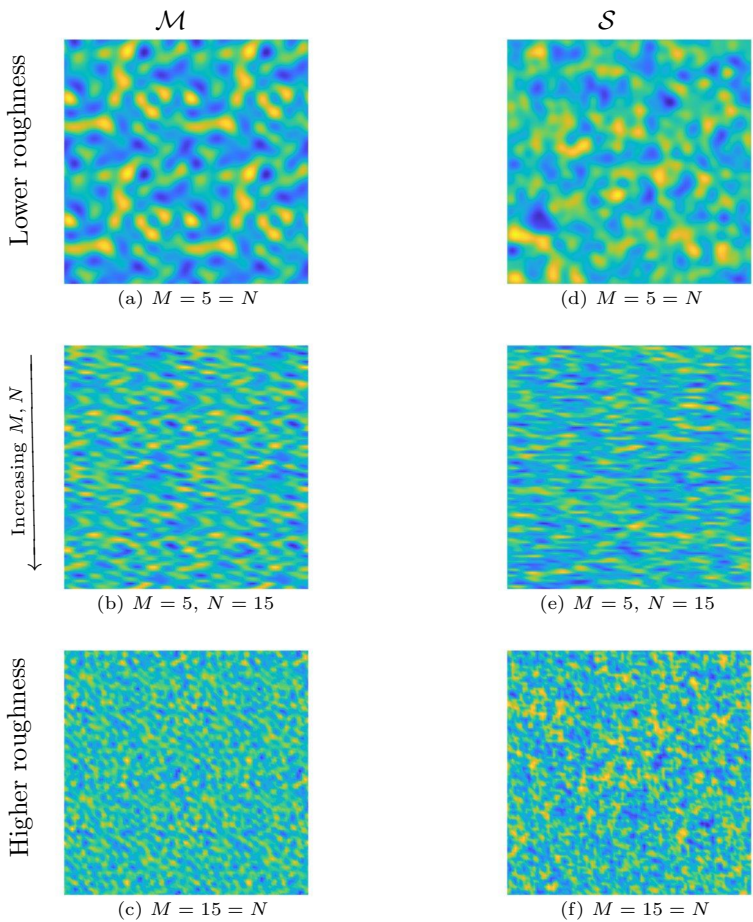

To determine \eref{eq:roughsurface1}, we first introduce $\tilde{a}_{m,n}\sim\mathcal{N}(0,1)$ and  $\phi_{m,n}\sim \calu(0,\pi)$.
This specifies the pre-surface function $\tilde{z}=\tilde{z}(x,y)$.
Next, we control the amplitude of the rough surface $\calm$ to be within a specific range $[-\MAmp,\MAmp]$ by scaling  according to
\begin{equation} \label{eq:roughsurface}
\begin{aligned}z:= \frac{\MAmp}{~\|\tilde{z}\|_{\infty}}\tilde{z} ,
\quad\mbox{i.e., } a_{m,n}:= \frac{\MAmp}{~\|\tilde{z}\|_{\infty}} \tilde{a}_{m,n}.
\end{aligned}
\end{equation}
Figures~\ref{fig:roughsurface}~(a)--(c), first column, show some rough surfaces $\M$ with amplitude $\MAmp=10^{-3}$ and various $(M,N)$.
The construction of the rough surfaces  $\mathcal{S}$ appearing in the second column will be discussed in Section \ref{sec:RRS}.

\textbf{Remark}: One can also add a decay condition with respect to frequencies $m,\,n$ and a frequency attenuation parameter $\beta$ to the pre-coefficient $\tilde{a}_{m,n}$ via
\[\tilde{a}_{m,n}\sim \frac{1}{({m}^2+n^2)^{\beta/2}} \,{\mathcal{N}(\mu,\sigma)}.\]

\subsection{Surface roughness and Laplace-Beltrami operator}

Applying the definition of  the Laplace-Beltrami operator {in} \eref{eq:LapBeltrami} to the rough surface function in \eref{eq:roughsurface1}--\eref{eq:roughsurface},  the relationship between the eigenvalues of the diffusion tensor and the geometry of the surface can be made explicit. Firstly, the {\emph{diffusion tensor} $\cala$} in \eref{eq:LapBeltrami} is
\begin{equation}\label{eq:D}
\cala(x,y) =
\begin{bmatrix}
\cala_1 & \cala_2 \\
\cala_2 & \cala_4
\end{bmatrix}(x,y)
:=\sqrt{g}G^{-1}=\frac{1}{\sqrt{g}}
\begin{bmatrix}
        1+z_y^2 & -z_x z_y \\
        -z_x z_y & 1+z_x^2 \\
\end{bmatrix},
\end{equation}
with the partial derivatives $z_x$ and $z_y$ computing from \eref{eq:roughsurface1}.
We are interested in obtaining the eigenvalues and eigenvectors of $\cala$.
{Instead of} computing these quantities directly, it turns out to be  easier to first compute the eigenvalues and eigenvectors of the Riemannian metric tensor $G$ in \eref{eq:G}.
From the characteristic equation of $G$
\[
|G-\lambda I|=(1+z_x^2-\lambda)(1+z_y^2-\lambda)-z_x^2z_y^2=0,
\]
we find that the eigenvalues $\lambda^G$ of $G$ are  $1+z_x^2+z_y^2=g$ and $1$.
Since
\begin{equation}
\label{eq:eig1}
    \cala \vec{v} = \sqrt{g}G^{-1} \vec{v}   = \frac{\sqrt{g}}{\lambda^G}\vec{v},
\end{equation}
any eigenvector $\vec{v}$ of $G$ is also an eigenvector of $\cala$ {
,} and the eigenvalues $\lambda^\cala$ of $\cala$ are
\begin{equation}
\label{eq:eig2}
    \lambda^{\cala}
    =\Big\{ \sqrt{g}=\sqrt{1+z_x^2+z_y^2} ,\; \frac{1}{\sqrt{g}} \Big\}.
\end{equation}
Note that, by \eref{eq:roughsurface}, the matrix function $[G-\lambda^G I](z_x,z_y) = C(\MAmp) [G-\lambda^G I](\tilde{z}_x,\tilde{z}_y)$ for some $\MAmp$-dependent  {constant} $C(\MAmp)$.  This implies that all eigendirections of $G$ (and $\cala$) are independent of the amplitude $\MAmp$.
Simple calculations show that
$\lambda^\cala=1+\mathcal{O}({M^2N^2\delta^2_{\M}})$ varies with $\MAmp$ nonlinearly by a bijection.  Thus, the contour lines (but not the height) of $\lambda^\cala$ are independent of $\MAmp$.

\begin{figure}
	\centering
	\setlength{\abovecaptionskip}{0pt}
	\setlength{\belowcaptionskip}{0pt}
	\begin{tabular}{cc}
	\begin{overpic}[width=0.4\textwidth]{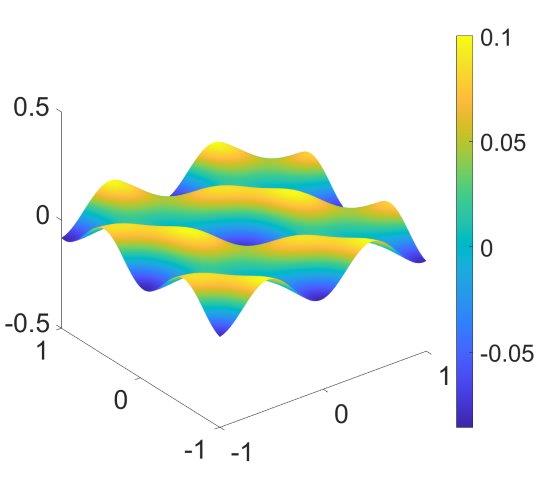}
		\put(18,83) {\scriptsize (a) $  M=1=N$  }
	\end{overpic}
	&
	\begin{overpic}[width=0.4\textwidth]{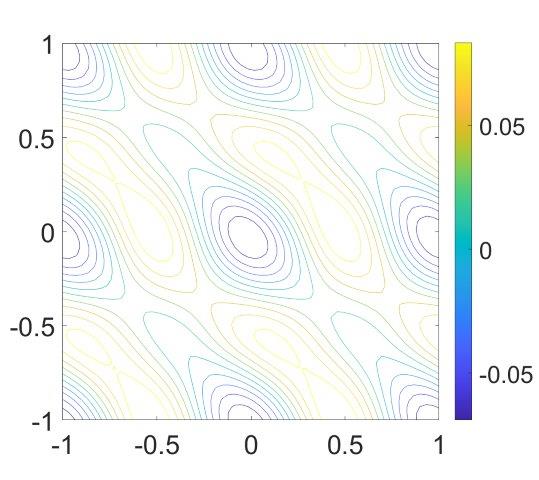}
		\put(10,83) {\scriptsize (b) Contour of a rough surface $\calm$}
	\end{overpic}
	\\
	\begin{overpic}[width=0.4\textwidth]{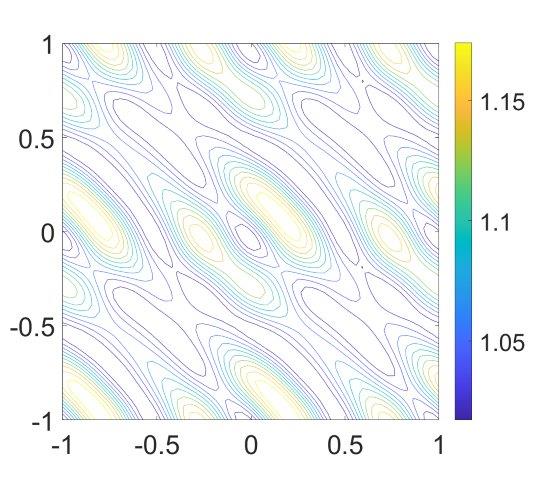}
		\put(18,83) {\scriptsize (c) maximum eigenvalue $\lambda^\cala_{\text{max}}$  }
	\end{overpic}
	&
	\begin{overpic}[width=0.4\textwidth]{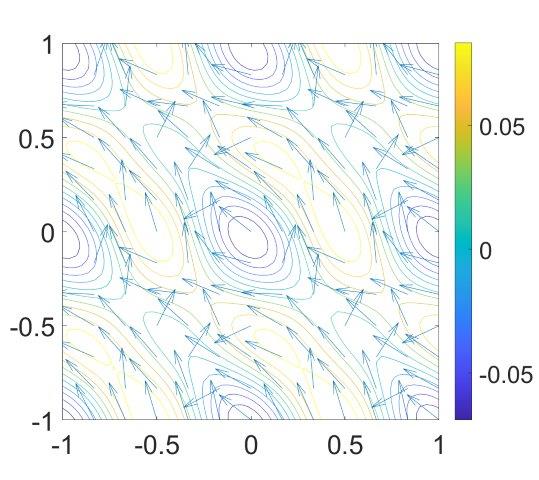}
		\put(18,83) {\scriptsize (d)  eigendirection of $\lambda^\cala_{\text{max}}$}
	\end{overpic}
	\\
		\begin{overpic}[width=0.4\textwidth]{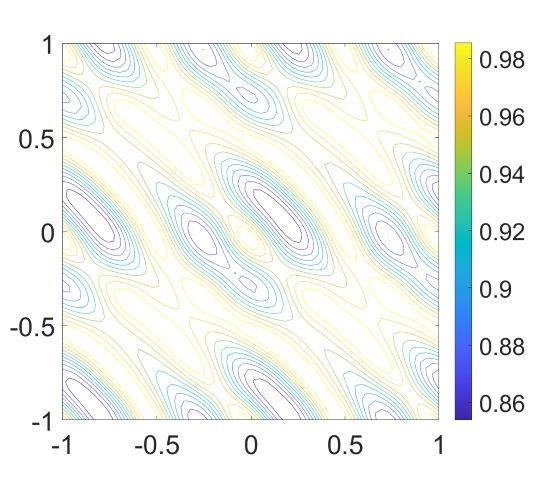}
		\put(18,83) {\scriptsize (e) minimum eigenvalue $\lambda^\cala_{\text{min}}$  }
	\end{overpic}
	&
	\begin{overpic}[width=0.4\textwidth]{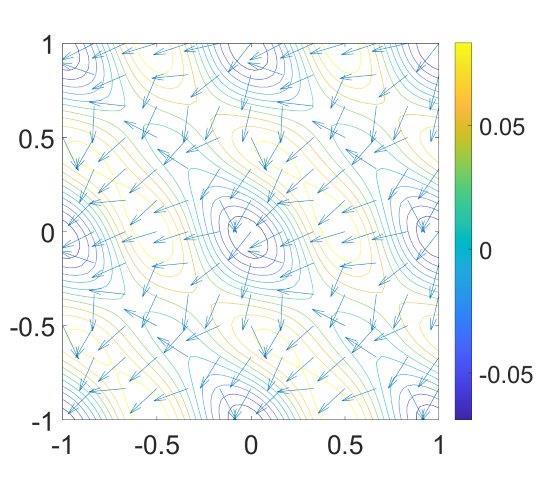}
		\put(18,83) {\scriptsize (f) eigendirection of $\lambda^\cala_{\text{min}}$ }
	\end{overpic}
	\end{tabular}
	\caption{ (a, b): 3D and contour plots of a rough surface {characterized by} \eref{eq:roughsurface}  with  $M=N=1$ and amplitude  $\MAmp=1E-1$; (c, d):  maximum eigenvalue and {corresponding} eigendirection of  $\cala$,  (e, f): minimum  eigenvalue and {corresponding} eigendirection of  $\cala$.}
	\label{fig:maxmineig}
\end{figure}

Figure~\ref{fig:maxmineig} illustrates the close relationship between the eigenvalues of $\cala$ and the geometric properties of a rough surface.
In particular, {when examining subfigures (b), (c), and (e),  it becomes evident that }the
maximum (and minimum) eigenvalues are larger (and smaller) at regions of the surface with the steepest  gradients.
In subfigures~(d) and (f), eigendirections are plotted on top of subfigure~(b).
We can clearly see that the eigenvectors corresponding to the maximum (and minimum) eigenvalues are tangent (and orthogonal) to the contour plots of eigenvalues. By adding roughness to the computational domain, an isotropic diffusion problem will become anisotropic and {directionality of gradients in the solution} will appear. Eigen-directions associated with large eigenvalues are preferred directions.

\section{Solving heat equations on rough surfaces $\calm$}\label{sec:FDM}

We begin by  constructing a finite difference scheme
for solving heat equations on rough surfaces $\calm \subset \R^3$ defined by \eref{eq:roughsurface1}--\eref{eq:roughsurface}. We work intrinsically by transforming the  PDEs on rough surfaces to the parameter space ${\cali^2=[-L,L]^2}\subset\R^2$.
Subsequently, we move on to solving reaction-diffusion systems for pattern formation.

We consider the heat equation on a rough surface {described by}
\begin{subequations}\label{eq heat}
\begin{equation}\label{eq:Parbolic}
\frac{\partial u}{\partial t} (\bx,t)   -\Delta_{\M}u(\bx,t)  = h(\bx,t)
\qquad  \text{for } \bx\in\M  \text{, } t\in(0,T],
\end{equation}
where $u:\M \times (0,T] \to \R$,
subject to  periodic boundary conditions  on $\partial {\M} := \partial \cali^2 \times z(\partial \cali^2)$  with $ z(\partial \cali^2)$ being the surface function in \eref{eq:roughsurface1}, i.e.,
 \begin{equation}\label{eq:HeatBC}
 \begin{aligned}
 u([-L,y,z(-L,y)]^T,t)&=u([L,y,z(L,y)]^T,t),
 \quad  \text{for }  y \in \cali , \ t\in(0,T],
\\
 u([x,-L,z(x,-L)]^T,t)&=u([x,L,z(x,L)]^T,t),
\quad \text{for }    x\in \cali , \ t\in(0,T].
 \end{aligned}
\end{equation}
\end{subequations}

Since all surface points are in the form $\bx =[x,y,z(x,y)]^T\in\M$ for $[x,y]^T\in \cali^2$ and $z$ in \eref{eq:roughsurface1}, we discretize the parameter space $\cali^2$ by some set of $n_X n_Y$ tensor-product grid points{,} $[X,Y]\in \R^{n_X n_Y \times 2}\subset \cali^2$. Let $u^j(x,y)\approx u( [x,y,z(x,y)]^T ,t^j)$ for   $[x,y]^T\in \cali^2$.
We begin by  temporal   discretization of  \eref{eq:Parbolic} by the $\theta$-method
\begin{eqnarray}\label{eq:timeBE}
\frac{u^{j+1}-u^{j}}{\tau}=\theta \big(\Delta_{\M}u^{j+1}+h ^{j+1}\big)+(1-\theta) \big(\Delta_{\M}u^{j}+h^{j}\big),
\end{eqnarray}
for an equispaced partition {$\{t^j\}_{j=0}^{M}$} of $[0,T]$ with {a} step-size $\tau$.
By some user-selected finite difference scheme for the first derivative, we construct differentiation matrices $\D_k \in \R^{n_X n_Y \times n_X n_Y }$ such that, for any $C^1$-function $w:\cali^2\to\R$ satisfying periodic boundary condition \eref{eq:HeatBC}, we have
\[
    \D_k w(X,Y) \approx \frac{\partial w}{\partial x_k}(X,Y), \qquad \text{for } k\in\{1,2\},
\]
where the $n_X n_Y \times 1$ vector $w(X,Y)$ (and $\frac{\partial w}{\partial x_k}(X,Y)$)  contains nodal values of $w$ (and $\frac{\partial w}{\partial x_k}$) at grid points in $[X,Y]$.
{By} working directly on {the} definitions \eref{eq:LapBeltrami}--\eref{eq:D}, we can approximate the Laplace-Beltrami term $\Delta_\M w$ at grids $[X,Y]$ by nodal function values $W :=w(X,Y)$ {as follows:}
\begin{equation}    \label{eq:disLB}
\begin{aligned}
\Delta_{\calm}w(X,Y) & \approx {\frac{1}{\sqrt{g}}(X,Y) \circledast} \bigg( \D_1  \Big(\cala_1(X,Y) \circledast  (\D_1 W )\Big)   +\D_2 \Big(  \cala_4(X,Y) \circledast (\D_2 W)\Big)
    \\ &\qquad +\D_1 \Big(\cala_2(X,Y) \circledast (\D_2 W  ) \Big)  +\D_2\Big( \cala_2(X,Y) \circledast( \D_1 W) \Big) \bigg)
    \\
& = {\frac{1}{\sqrt{g}}(X,Y) \circledcirc} \Big[ \D_1  \Big(\cala_1(X,Y) \circledcirc  \D_1\Big)
         +\D_2 \Big(  \cala_4(X,Y) \circledcirc \D_2\Big)
    \\ &\qquad\qquad +\D_1 \Big(\cala_2(X,Y) \circledcirc \D_2  \Big)
            +\D_2\Big( \cala_2(X,Y) \circledcirc \D_1 \Big) \Big] W
    \\
& =: \triangle_{\calm,h} W,
    \end{aligned}
\end{equation}
where $\triangle_{\calm,h} \in \R^{n_X n_Y \times n_X n_Y }$,  $\circledast$ 
is the element-wise Hadamard product of matrices and $\circledcirc$ is a vector-matrix product defined {as} $\vec a \circledcirc [\vec b_1, \cdots,   \vec b_n]:= [\vec a\circledast \vec b_1, ..., \vec a\circledast \vec b_n]$. See {the} Appendix for the detailed construction of differentiation matrices $\D_1,\D_2$. Combining \eref{eq:timeBE} and \eref{eq:disLB} yields a fully discretized scheme for the update in time.

\begin{figure}
	\centering
	\setlength{\abovecaptionskip}{0pt}
	\setlength{\belowcaptionskip}{0pt}
	\begin{tabular}{cc}
	\begin{overpic}[width=0.30\textwidth]{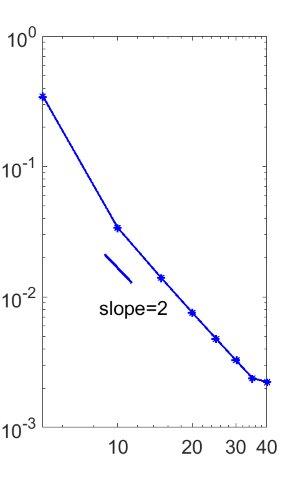}
	\put(30,95) { (a) }
	\put(27,3) { $n_X$}
\put(-5,50){\rotatebox{91}{ $\ell^{2}$ error}}
	\end{overpic}
	&
	\begin{overpic}[width=0.30\textwidth]{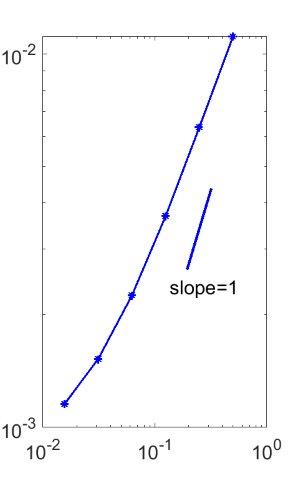}
	\put(30,95) { (b) }
	\put(33,3) { $\tau$}
	\end{overpic}
	\end{tabular}
	\caption{Accuracy and convergence results for the heat equation on the rough surface in Figure~\ref{fig:maxmineig}(a) with $M=N=1, \delta_M=1E-2$: (a) Convergence with respect to spatial refinement. (b) Convergence with respect to time.}
	\label{fig:ex30}
\end{figure}
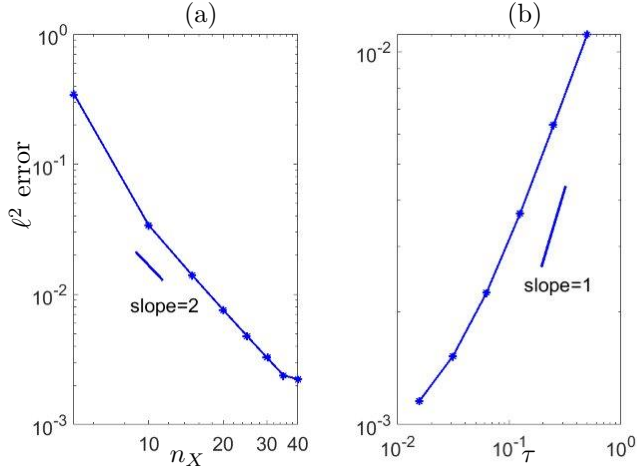

{To verify the accuracy and convergence of the introduced finite difference scheme, we consider a heat equation with a known exact solution.} The parameter space is set to be $\mathcal{I}^2=[-1,1]^2$. The exact solution is $u^*(x,y,t)=\exp(t)\sin(\pi x)\sin(\pi y)$. The function $f(x,y,t)$ and initial condition are generated from the exact solution $u^*$.  In all cases, we use the backward Euler method. Periodic boundary conditions are imposed.
The relative $\ell^{2}$-error at $t^n$ is computed as
\begin{equation}\label{ell2err}
\ell^{2}\text{-error} =\frac{\| u_N^n-u^* \|_2}{ \| u^* \|_2} .
\end{equation}
For the surface $M=N=1$ in Figure \ref{fig:maxmineig}~(a) with  $\delta_M=1E-2$,  Figure~\ref{fig:ex30}~(a) exhibits second order spatial convergence for ${n_X,n_Y\in\{5,10,\ldots,40\}},\ \tau=1E-3,\ T=1E-1$. Figure~\ref{fig:ex30}~(b) exhibits first order convergence  with respect to time for $n_X=90=n_Y,\ {\tau\in\{1/2,1/2^2,\ldots, 1/2^6\}},\ T=1$.

\subsection{Visualizing heat flow on a rough surface $\calm$}
\begin{figure}
	\centering
	\setlength{\abovecaptionskip}{0pt}
	\setlength{\belowcaptionskip}{0pt}
	\begin{tabular}{cc}
	\begin{overpic}[width=0.42\textwidth]{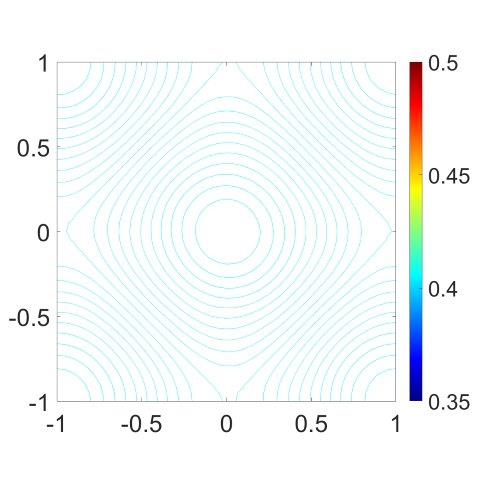}
	\put(40,90) { (a)    }
	\put(-5,30){\rotatebox{91}{ $\delta_\calm=0.1$}}
	\end{overpic}
	&
	\begin{overpic}[width=0.42\textwidth]{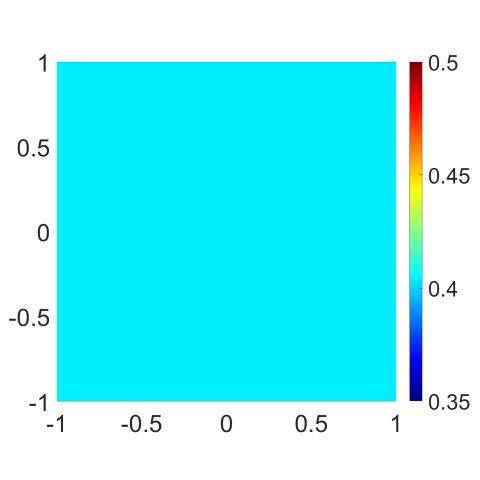}
	\put(40,90) { (b)    }
	\end{overpic}
	\\
		\begin{overpic}[width=0.42\textwidth]{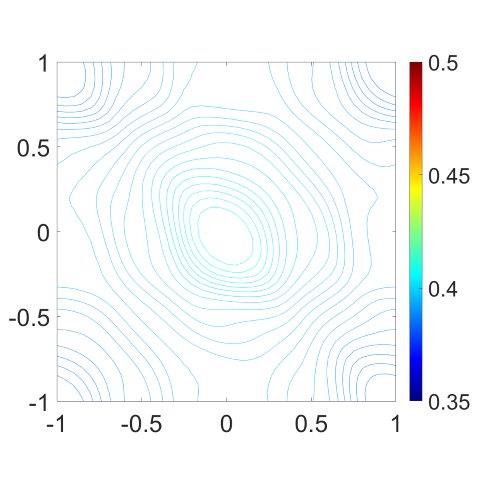}
		\put(40,90) { (c)  }
	\put(-5,30){\rotatebox{91}{ $\delta_\calm=0.5$}}
	\end{overpic}
	&
	\begin{overpic}[width=0.42\textwidth]{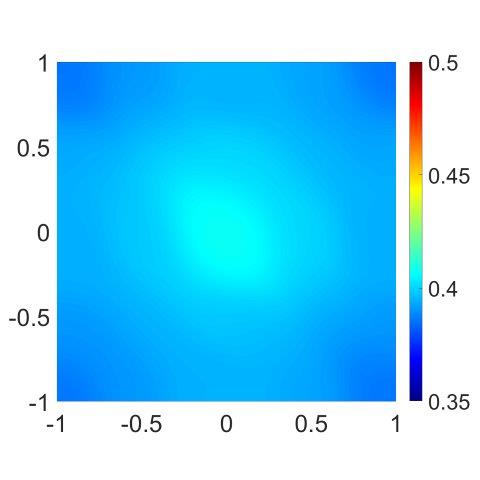}
	\put(40,90) { (d)}
	\end{overpic}
\\
	\begin{overpic}[width=0.42\textwidth]{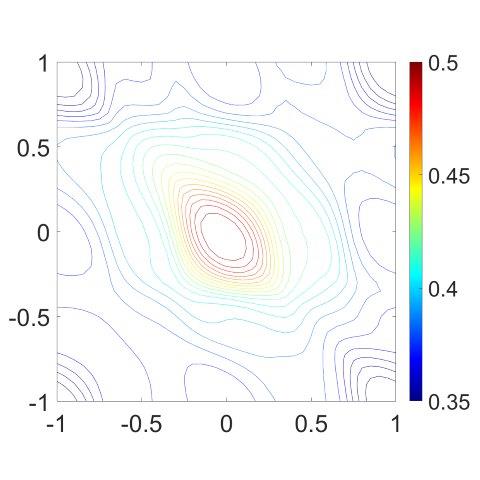}
	\put(40,90) { (e)}
	\put(-5,30){\rotatebox{91}{ $\delta_\calm=1$}}
	\end{overpic}
	&
	\begin{overpic}[width=0.42\textwidth]{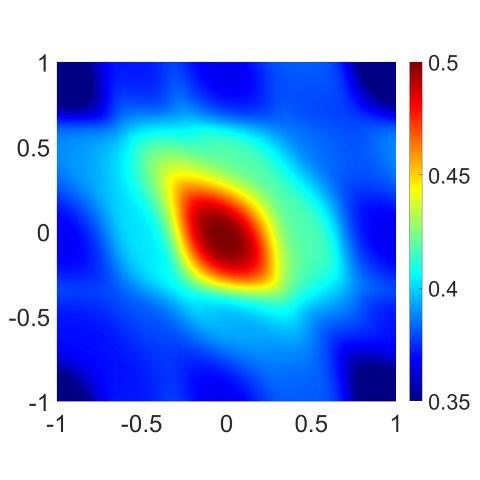}
		\put(40,90) { (f)}
	\end{overpic}
	\end{tabular}
	\caption{ Numerical {investigation of heat equation} {solutions on $\mathcal{M}$ with varying amplitudes $\delta_\mathcal{M}$}: (a, b): $\delta_\mathcal{M}=0.1$, (c, d): $\delta_\mathcal{M}=0.5$, (e, f): $\delta_\mathcal{M}=1$.  In all cases, $\tau=1E-3,\ T=1, \ n_X=41=n_Y$.}
	\label{fig:numheat}
\end{figure}
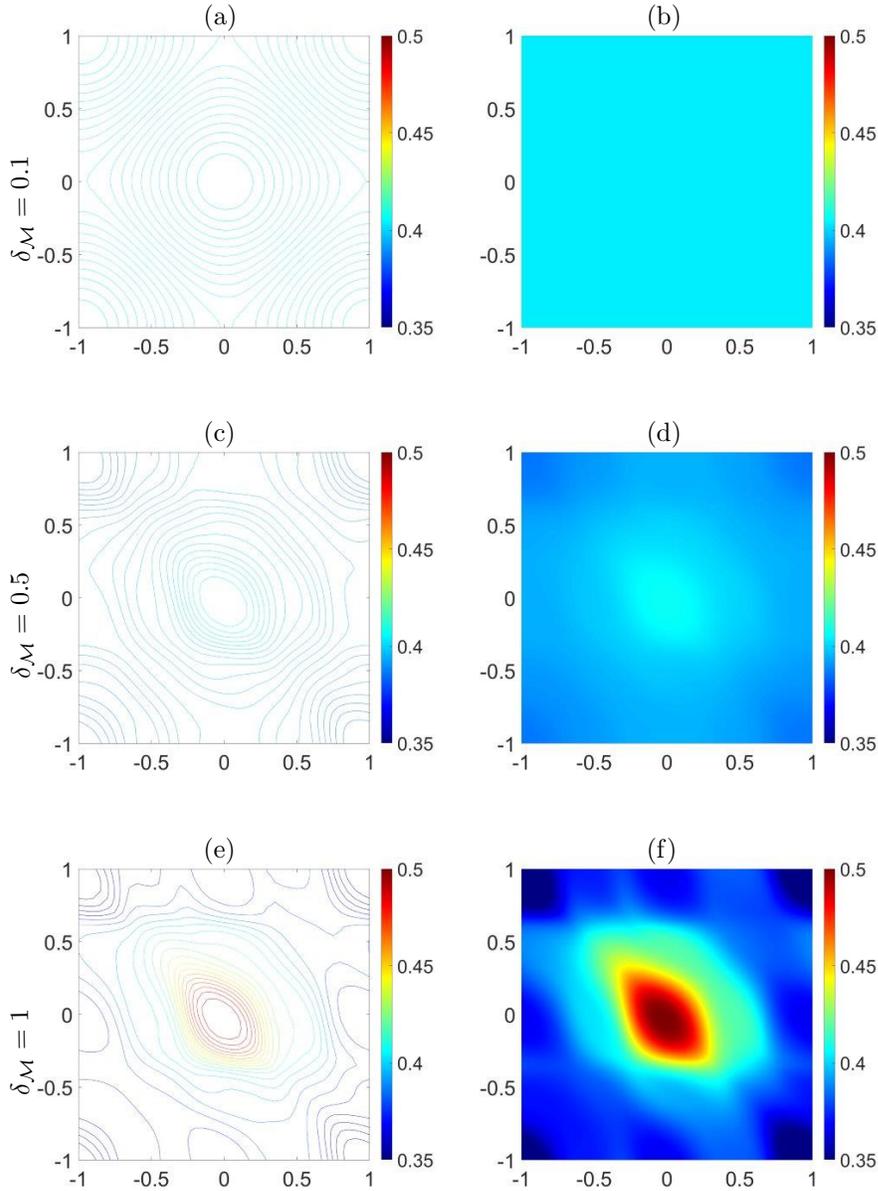

We {visualize} the heat flow under different amplitudes of rough surface $\calm$ in {Figure \ref{fig:maxmineig}} by solving the heat equations \eref{eq:Parbolic} with zero flux $h(\xi,t)=0$ and periodic boundary conditions \eref{eq:HeatBC}. The rough surfaces $\calm$ are defined over the parameter space  ${\cali^2=[-1,1]^2}\subset\R^2$.
The compatible initial condition is given by
\begin{equation} \label{eq heat ic}
u(x,y,z,0) = \cos(\pi x/2)\cos(\pi y/2).
\end{equation}
Figure~\ref{fig:numheat} shows the numerical solutions on the surface $\mathcal{M}$ in Figure \ref{fig:maxmineig}~(b) with $\delta_\mathcal{M}\in\{0.1,0.5,1\}$ in the respective rows. In all {cases}, we
set $\tau=1E-3,\ T=1$ and $n_X=41=n_Y$.
Figure~\ref{fig:numheat} shows that the
range of the numerical solutions becomes larger as the amplitude of $\mathcal{M}$ increases from $\delta_\mathcal{M}=0.1$ to $\delta_\mathcal{M}=1$. From the rough surface in
Figure \ref{fig:maxmineig}~(b) and solutions in Figure~\ref{fig:numheat},  it can be concluded that the heat flow, when projected {onto} the plane, is greatest in flat regions.

\section{{Pattern formation} on a rough surface $\calm$}\label{sec:num}
In this section, we {examine}
reaction-diffusion systems (RDS)  on a rough surface $\M$, {given by}
\begin{equation}
\left\{
\begin{aligned}
\partial_t u&=\delta_u\Delta_{\M}u +f_u(u,v),\\
\partial_t v&=\delta_v\Delta_{\M}v +f_v(u,v),
\end{aligned}
\right.
\label{eq:RDS}
\end{equation}
for some (concentration) functions $u,v: {\calm}\times(0,T]\to\R$, and reaction terms
\begin{equation}
\left\{
\begin{aligned}
f_u(u,v)&=\alpha u(1-\xi_1 v^2)+v(1-\xi_2 u),\\
f_v(u,v)&=\beta v\Big(1+\frac{\alpha \xi_1}{\beta} u v\Big)+u(\gamma+\xi_2 v),
\label{eq:RDSmodel}
\end{aligned}
\right.
\end{equation}
with parameters $\delta_u,\delta_v,\alpha, \beta, \xi_1, \xi_2$.
Pattern formation on very smooth (and usually closed) surfaces with little variation {has been} discussed in {\cite{Cheung+Ling-Kernembemethconv:18,shankar2020robust,lehto2017radial,mcdonald2007global}}{. However, to the best of the authors' knowledge, the study of pattern formation on rough surfaces has not been undertaken .}
In this section, we  aim to analyze the effect of  rough surfaces on pattern formation generated by the reaction-diffusion system~(\ref{eq:RDS})-(\ref{eq:RDSmodel}). {We are not primarily focused on numerical methods, and }simply extend the finite difference scheme {from} the previous section to work here. Other options for {discretizing} reaction-diffusion systems include the
finite element method  \cite{Guanghui2015Moving,2008Numerical,2013Implicit}, and various types of meshfree methods  \cite{AAMM-12-1327,2021ASLI}.

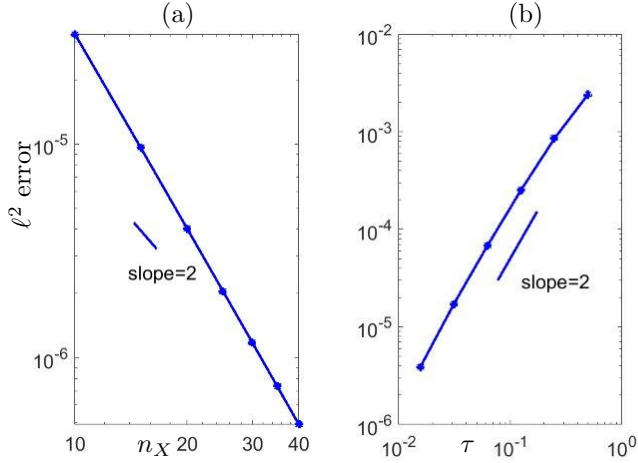
\begin{figure}
	\centering
	\setlength{\abovecaptionskip}{0pt}
	\setlength{\belowcaptionskip}{0pt}
	\begin{tabular}{cc}
	\begin{overpic}[width=0.30\textwidth]{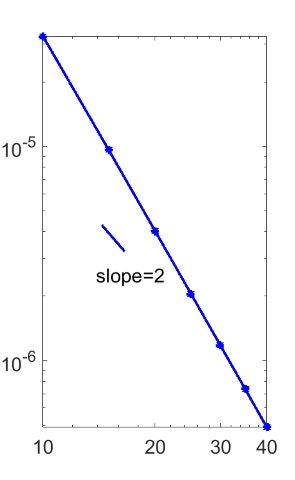}
	\put(25,95) { (a) }
	\put(20,5) { $n_X$}
    \put(-5,50){\rotatebox{90}{ $ \ell^{2}$ error}}
	\end{overpic}
	&
	\begin{overpic}[width=0.30\textwidth]{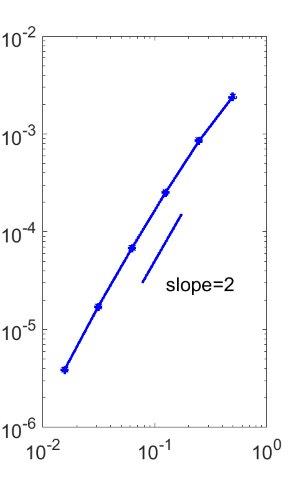}
	\put(30,95) { (b) }
	\put(20,5) { $\tau$}
	\end{overpic}
	\end{tabular}
	\caption{Accuracy and convergence results for a reaction-diffusion system on the rough surface  in Figure~\ref{fig:maxmineig}(a) with $M=N=1,\delta_M=1E-2$.
(a) Convergence with respect to spatial refinement. (b) Convergence with respect to time.
}
	\label{fig:ex40}
\end{figure}

As in Section~\ref{sec:FDM},
we work on an equispaced temporal partition $\{t^j\}_{j=0}^{N_T}$ with some {time step} $\tau>0$ and  tensor-product grid points $[X,Y]\in \R^{n_X n_Y \times 2}\subset \cali^2$.
Let
\[
    U^j\approx u( X,Y,z(X,Y), t^j) \quad \text{and} \quad V^j\approx v( X,Y,z(X,Y), t^j), \qquad 1\leq j \leq N_T,
\]
be the unknown nodal values we seek based on initial conditions $U^0$ and $V^0$.
Using the second order backward differentiation formula (BDF2) \cite{2009Computer,1994Implicit} and \eref{eq:disLB}  to discretize the RDS \eref{eq:RDS} with periodic boundary condition \eref{eq:HeatBC},  we obtain the following fully-discrete system of equations on $\calm$
\begin{equation}\label{eq:RDSSEMIDIS}
\left\{
\begin{split}
& 3U^{j+1}- 2 \tau \delta_u \triangle_{\calm,h} U^{j+1} =
4\tau f_u\Big(U^{j},V^{j}\Big)-2\tau f_u\Big(U^{j-1},V^{j-1}\Big)+4U^j-U^{j-1},
\\
&  3V^{j+1}- 2 \tau \delta_v  \triangle_{\calm,h} V^{j+1} =
4\tau f_v\Big(U^{j},V^{j}\Big)-2\tau f_v\Big(U^{j-1},V^{j-1}\Big)+4V^j-V^{j-1},
\end{split}
\right. 
\end{equation}
for $1 \leq j \leq N_T$,
subject to some  yet-to-be specified (usually random) initial condition and some first order approximations to the solutions at the first time step, $U^1$ and $V^1$.
Note that the two equations in \eref{eq:RDSSEMIDIS} are not coupled,  { and} the computational cost is of the same order as that of solving two scalar heat equations.
Here, we compute {the} RDS with exact solutions to {demonstrate} the convergence behavior of the finite difference scheme. The parametric domain is chosen to be $\mathcal{I}^2=[-1,1]^2$. The exact solutions are
\[u^*(x,y,t)=\exp(t)\sin(2\pi x)\sin(\pi y),\quad v^*(x,y,t)=\exp(t)\sin(\pi x)\sin(2\pi y).\]
In order to obtain RDS with these exact solutions, we construct  new reaction functions $F_u(u,v)$ and $F_v(u,v)$ as
\[F_u(u,v)=f_u(x,y,t)+f_u(u,v),\quad F_v(u,v)=f_v(x,y,t)+f_v(u,v),\]
with
\begin{equation*}
\begin{aligned}
f_u(x,y,,t)&=\frac{\partial u^*(x,y,,t)}{\partial t}-\delta_u \triangle_{\mathcal{M}} u^*(x,y,t)-f_u(u^*,v^*), \\
 f_v(x,y,t)&=\frac{\partial v^*(x,y,t)}{\partial t}-\delta_v \triangle_{\mathcal{M}} v^*(x,y,t)-f_v(u^*,v^*).
 \end{aligned}
\end{equation*}
Periodic boundary conditions are imposed and initial conditions are generated from the exact solution $u^*$, $v^*$. The relative $\ell^{2}$-errors at $t^n$ are computed by those of $u$ and $v$ as in \eref{ell2err} according to
\[ \ell^{2}\text{-error} =\frac{\sqrt{(\ell^{2}\text{-error of $u$})^2+ (\ell^{2}\text{-error of $v$})^2}}{n_E} .\]
We consider the surface in Figure \ref{fig:maxmineig}~(a) with  $M=N=1$ and $\delta_M=1E-2$. In Figure~\ref{fig:ex40}~(a), we {observe} second order spatial convergence for ${n_X,n_Y\in\{10,15,\ldots,40\}}$, $ \tau=1E-3$, $T=1E-1${.}  Figure~\ref{fig:ex40}~(b) shows second order convergence in time for $n_X=90=n_Y,\ {\tau\in\{1/2,1/2^2,\ldots, 1/2^6\}},\ T=1$.

\begin{table}
	\centering
	\caption{ Parameters of the reaction-diffusion system \eref{eq:RDS}-\eref{eq:RDSmodel} for generating spots and stripes patterns on rough surfaces. }
	\begin{tabular}{|c|c|c|c|c|c|c|c|}
	\hline
	Pattern &  $\delta_v$ &  $\delta_u$& $\alpha$ & $\beta$ & $\gamma$& $\xi_1$& $\xi_2$  \\
	\hline
	Spots   &  $10^{-3}$ & $0.516 \delta_v$ & $0.899$  & $-0.91$ & $-0.899$ & $0.02$ & $0.2$  \\
\hline
	Stripes & $10^{-3}$  & $0.516 \delta_v$ & $0.899$ & $-0.91$ & $-0.899$ & $3.5$ & $0$ \\
	\hline
	\end{tabular}
	\label{tbl:table1}
\end{table}

\subsection{The pattern generation process}
In this subsection, we show the formation of irregular patterns as we {transition} from a flat two dimensional domain to rough surfaces with {varying} amplitudes.  We {begin} with patterns on the flat domain $[-1,1]^2$ as in \cite{AAMM-12-1327}, i.e., the rough surface with zero amplitude  ($\MAmp=0$). {The model and} surface parameters are {determined based on the values in Table~\ref{tbl:table1}.}
{For discretization,} we select grid parameters $n_X=90,\ n_Y=n_X$, and a time step-size $\tau=0.5$.  The rough surfaces $\M$ with $M=N=5$ are used in this part.

\begin{figure}
	\centering
	\setlength{\abovecaptionskip}{0pt}
	\setlength{\belowcaptionskip}{0pt}
\begin{tabular}{cc} \begin{overpic}[width=0.37\textwidth]{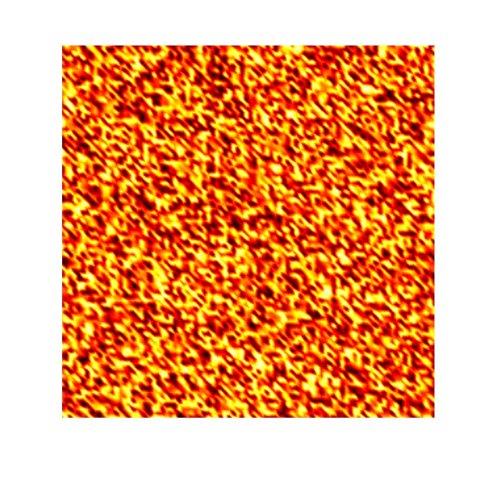}
	\put(40,95) {Spots}
	\put(0,30){\rotatebox{91}{$\MAmp= 0,\ t=0$}}
	\end{overpic}
	&	
	\begin{overpic}[width=0.37\textwidth]{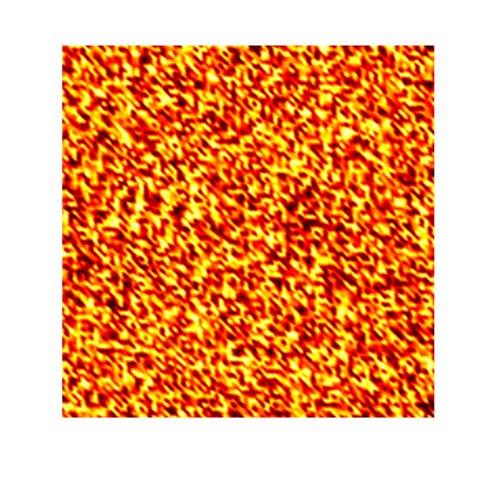}
	\put(40,95) {Stripes}
	\put(0,30){\rotatebox{91}{$\MAmp= 0,\ t=0$}}
	\end{overpic}
\\
\begin{overpic}[width=0.37\textwidth]{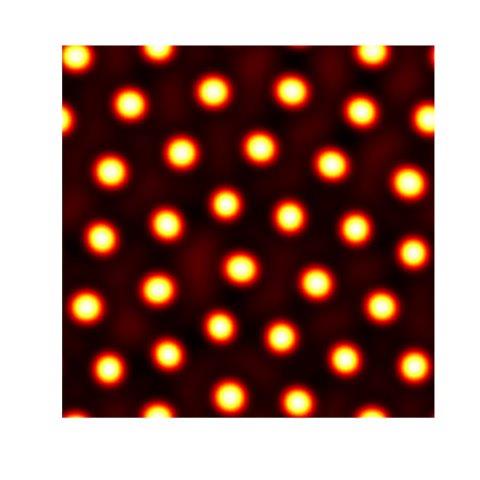}
	\put(0,30){\rotatebox{91}{$\MAmp= 0,\ t=800$}}
	\end{overpic}
&
\begin{overpic}[width=0.37\textwidth]{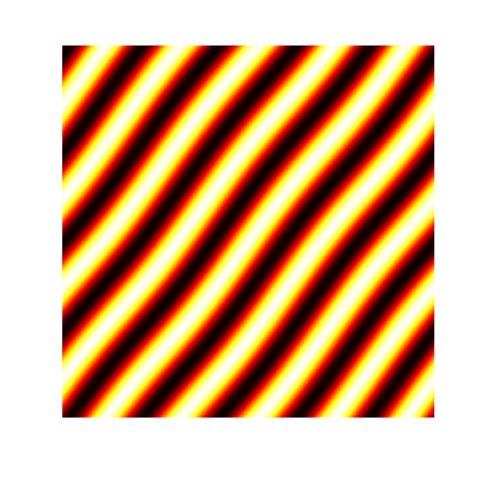}
	\put(0,30){\rotatebox{91}{$\MAmp= 0,\ t=4000$}}
	\end{overpic}
\\
\multicolumn{2}{c}{
{\rotatebox{-91}{ $\xrightarrow{\quad \rotatebox{91}{\text{ Increase
$\delta_\calm$ by   0.01 at a time
}}\quad}$ }}
}
\\
\begin{overpic}[width=0.37\textwidth]{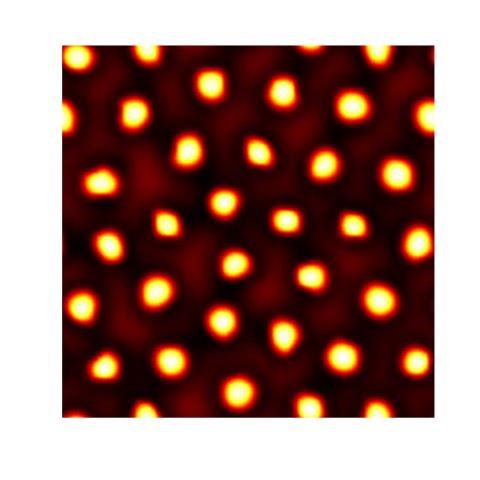}
	\put(0,12){\rotatebox{91}{$\MAmp= 0.1 ,\ t=800\times 11$}}
    \end{overpic}
&
\begin{overpic}[width=0.37\textwidth]{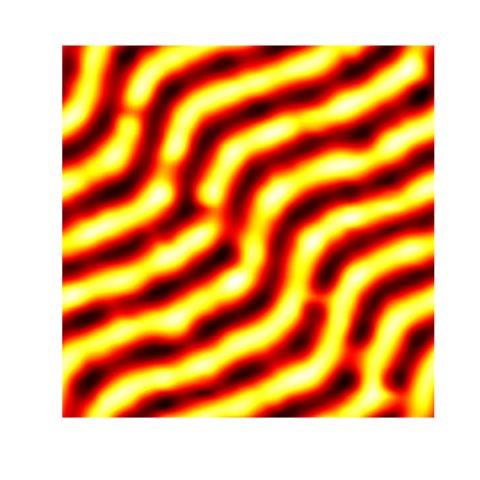}
	\put(0,12){\rotatebox{91}{$\MAmp= 0.1 ,\ t=4000\times 11$}}
	\end{overpic}
\end{tabular}
	\caption{Pattern generation on rough surfaces $\M$ with $M=N=5$ and amplitude $\MAmp$ increasing from $0$ to $0.1$. Parameters for spots and stripes are set according to Table \ref{tbl:table1}. {The discretization sets as} $n_X=90=n_Y, \ \tau=0.5$.
}\label{fig:ex31}
\end{figure}

The first row of Figure~\ref{fig:ex31} plots the initial conditions used to compute the spots and stripes patterns, respectively.  These are random values generated within the interval $[-0.5,0.5]$. The second {row} of Figure~\ref{fig:ex31} shows the steady state patterns on the surface with zero amplitude. Perfect spots and stripes are obtained{, similar} to our previous results in \cite{AAMM-12-1327}.
Next, we set the initial conditions for the next amplitude
to be the steady solutions from the zero amplitude rough surface. {In other words,} we use the solution of spots with $\MAmp=0,\ T=800$ and the solution of stripes with $\MAmp=0,\  T=4000$. By increasing the amplitude of the rough surface from $\MAmp=0$ to $\MAmp=0.1$ with increments of $0.01$, and setting initial conditions using the previous steady state, we achieve the final patterns for $\MAmp=0.1$.  These are shown in the third row of Figure~\ref{fig:ex31}. For rough surfaces under {a} small amplitude $\MAmp=0.05$, the spots and stripes are similar to those with $\MAmp=0$. However, for {a} larger amplitude $\MAmp=0.1$, both spots and stripes become irregular. {From this,} we can conclude that the steady state patterns become irregular as the amplitude of the rough surface $\M$ increases.

\subsection{Patterns on surfaces with different spatial frequencies and amplitudes} Properties of rough surfaces are influenced by the amplitude $\MAmp$ and the spatial frequencies  $M,\ N$ in \eref{eq:roughsurface}.
To better understand patterns on rough surfaces with different properties, we compute steady state patterns for  \[\MAmp\in\{0.05,0.1\}, \quad (M,N)\in\{ (5,15),\ (15,15)\},\]
on the parameter space ${\cali^2=[-0.5,0.5]^2} \subset\R^2$. Periodic boundary conditions on $\partial {\M} := \partial \cali^2 \times z(\partial \cali^2)$ are prescribed.

\begin{figure}
	\centering
	\setlength{\abovecaptionskip}{0pt}
	\setlength{\belowcaptionskip}{0pt}
	\begin{tabular}{cc}	
\begin{overpic}[width=0.45\textwidth]{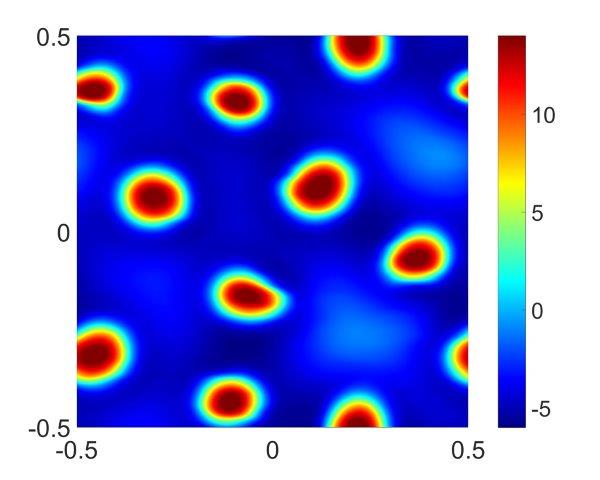}
			\put(38,78) { 2-D view}
	\put(0,20){\rotatebox{91}{ $M=5,\,N=15$}}
	\end{overpic}
	&
	\begin{overpic}[width=0.45\textwidth]{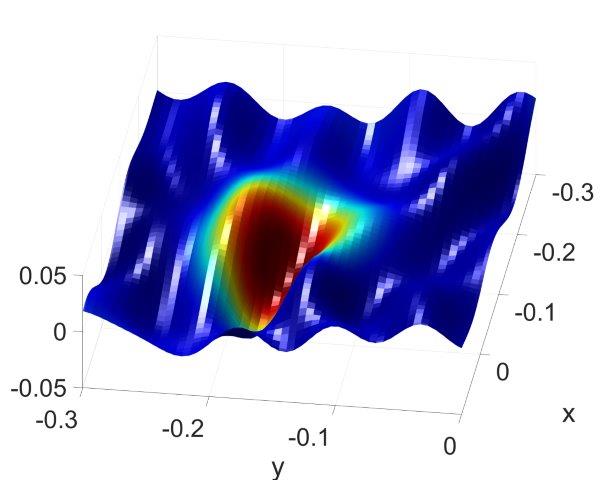}
    			\put(20,78) { Zoom in 3-D view}
		\end{overpic}
		\\
\begin{overpic}[width=0.45\textwidth]{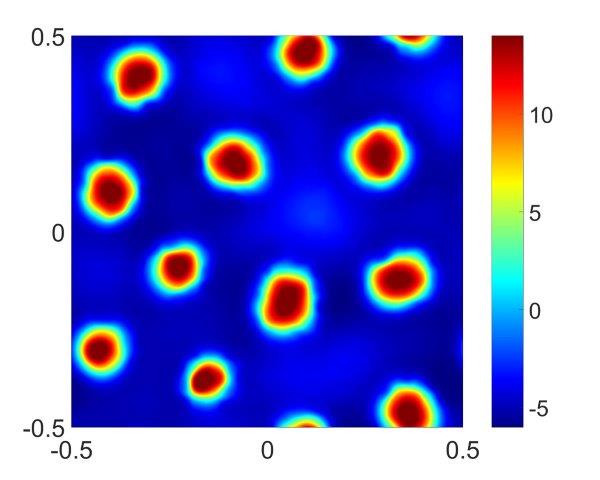}
	\put(0,20){\rotatebox{91}{ $M=15=N$}}
	\end{overpic}
	&
	\begin{overpic}[width=0.45\textwidth]{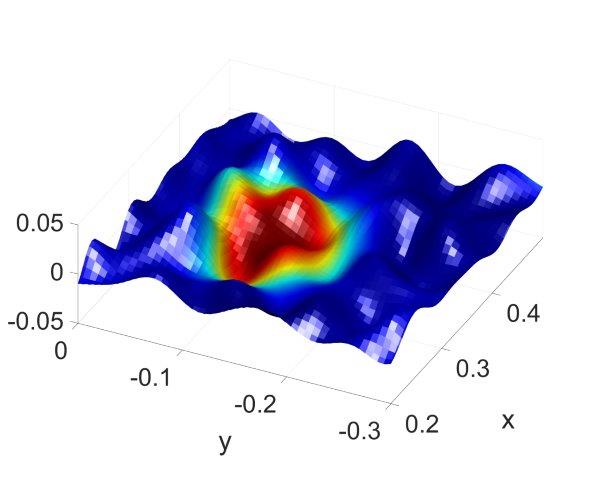}
		\end{overpic}
	\end{tabular}
	\caption{Pattern {formation} on rough surfaces with {an} amplitude $\MAmp=0.05$. Discretization parameters are $n_X=170=n_Y,\ \tau=0.5,\ T=800, \ \cali^2:=[-0.5,0.5]^2$.}
	\label{fig:ex321}
	\end{figure}

\begin{figure}
	\setlength{\abovecaptionskip}{0pt}
	\setlength{\belowcaptionskip}{0pt}
	\begin{tabular}{cc}	
	\centering
\begin{overpic}[width=0.45\textwidth]{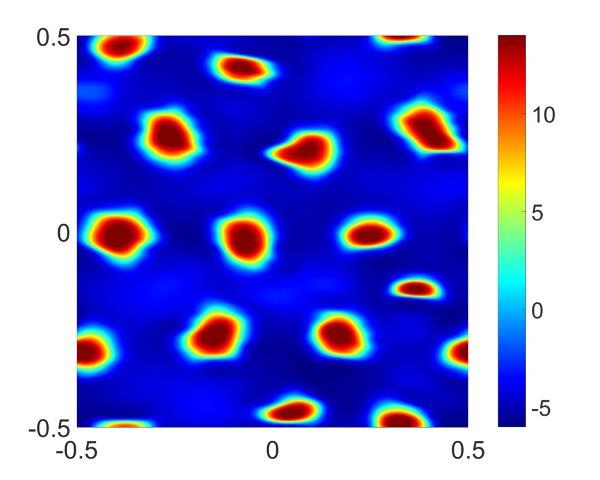}
		\put(38,78) { 2-D view}
	\put(0,20){\rotatebox{91}{ $M=5,\,N=15$}}
	\end{overpic}
	&	
	\begin{overpic}[width=0.45\textwidth]{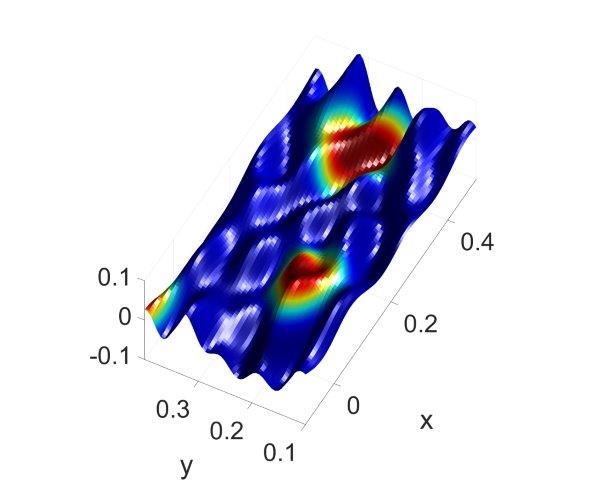}
    			\put(20,78) { Zoom in 3-D view}
		\end{overpic}
\\
\begin{overpic}[width=0.45\textwidth]{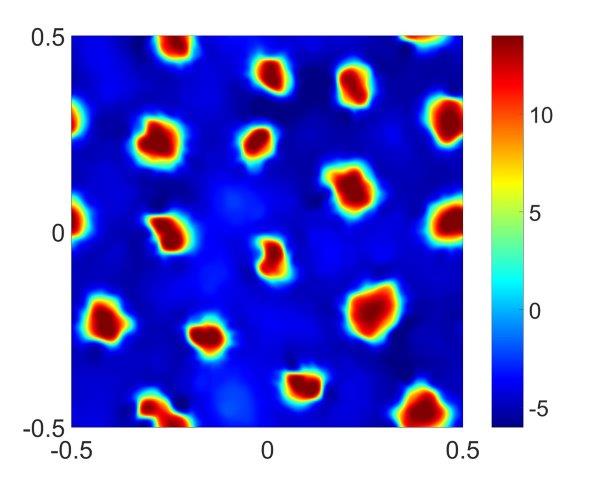}
		\put(0,20){\rotatebox{91}{ $M=15=N$}}
	\end{overpic}
	&	\begin{overpic}[width=0.45\textwidth]{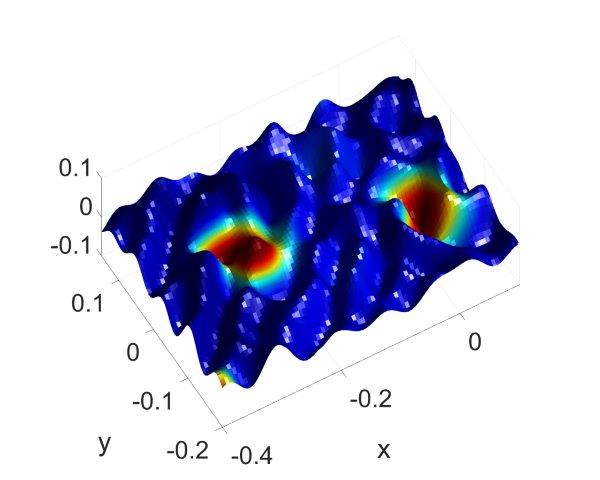}
		\end{overpic}
	\end{tabular}
	\caption{Patterns  on rough surfaces {$\mathcal{M}$ with an} amplitude $\MAmp =0.1$. Discretization parameters are $n_X=170=n_Y,\tau=0.5, T=800, \ \cali^2:=[-0.5,0.5]^2$.}
	\label{fig:ex322}
	\end{figure}

\begin{figure}
	\centering
	\setlength{\abovecaptionskip}{0pt}
	\setlength{\belowcaptionskip}{0pt}
	\begin{tabular}{cc}
\begin{overpic}[width=0.45\textwidth]{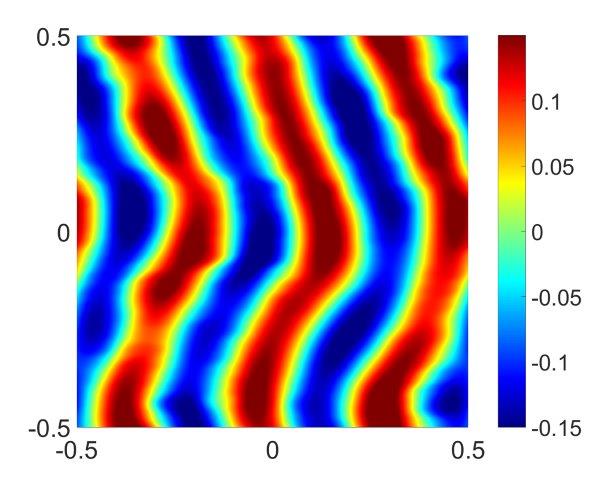}
			\put(38,78) { 2-D view}
	\put(0,20){\rotatebox{91}{ $M=5,\,N=15$}}
	\end{overpic}
	&
	\begin{overpic}[width=0.45\textwidth]{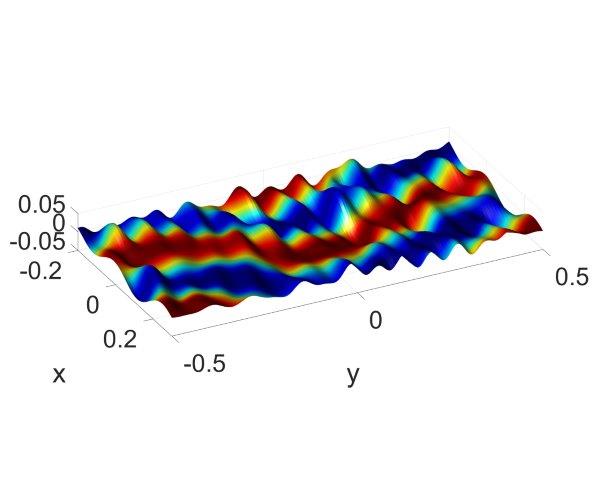}
    		\put(20,78) { Zoom in 3-D view}
		\end{overpic}
		\\
\begin{overpic}[width=0.45\textwidth]{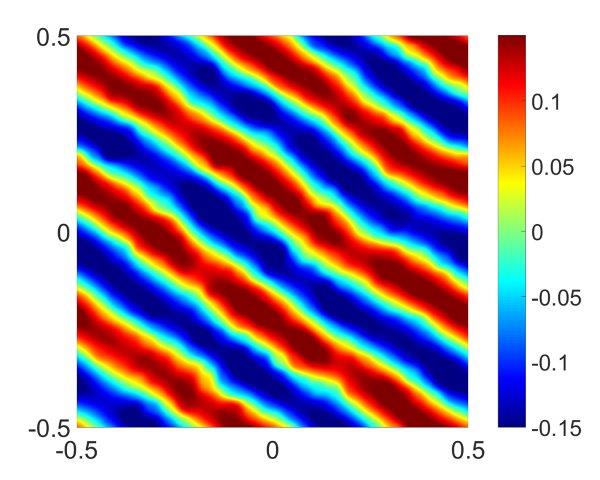}
	\put(0,20){\rotatebox{91}{ $M=15=N$}}
	\end{overpic}
	&
	\begin{overpic}[width=0.45\textwidth]{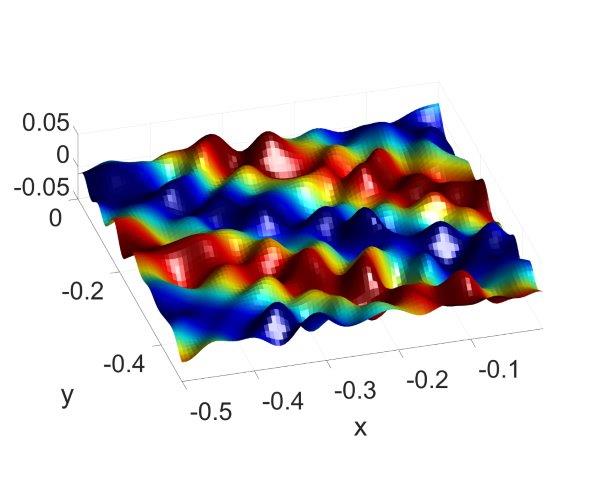}
		\end{overpic}
	\end{tabular}
	\caption{ Patterns on rough surfaces {$\mathcal{M}$ with an} amplitude $\MAmp \in\{0.05\}$ and $(M, N)\in \{(5, 15), (15,15)\}$. The model parameters for stripes are set according to Table~\ref{tbl:table1}. Discretization parameters are $\ n_X=170=n_Y,  \tau=0.5,\ T=4000$.}
	\label{fig:ex331}
	\end{figure}

\begin{figure}
	\centering
	\setlength{\abovecaptionskip}{0pt}
	\setlength{\belowcaptionskip}{0pt}
	\begin{tabular}{cc}
\begin{overpic}[width=0.45\textwidth]{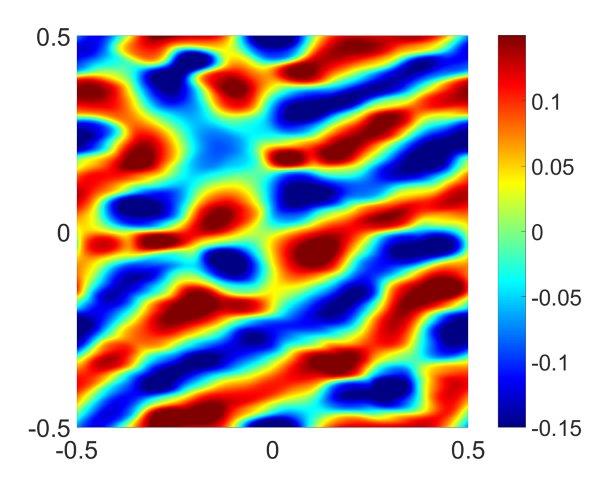}
			\put(38,78) { 2-D view}
	\put(0,20){\rotatebox{91}{ $M=5,\,N=15$}}
	\end{overpic}
	&
	\begin{overpic}[width=0.45\textwidth]{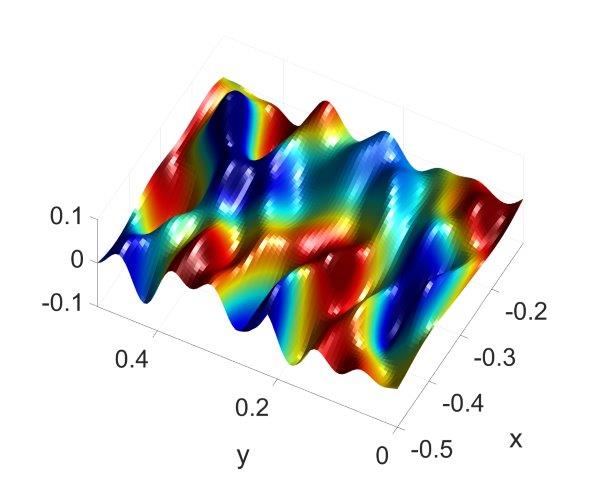}
    		\put(20,78) { Zoom in 3-D view}
		\end{overpic}
		\\
\begin{overpic}[width=0.45\textwidth]{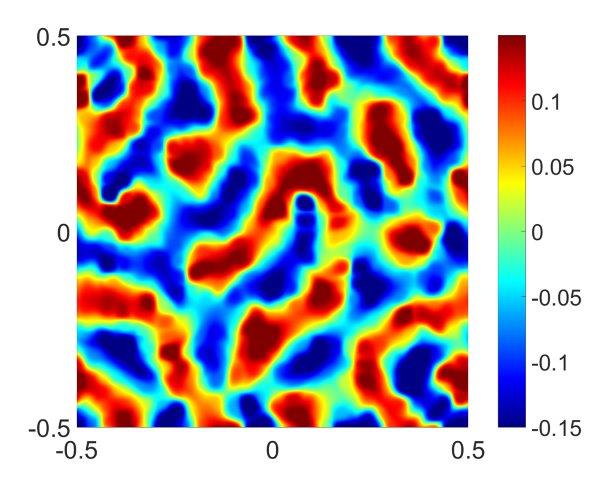}
	\put(0,20){\rotatebox{91}{ $M=15=N$}}
	\end{overpic}
	&
	\begin{overpic}[width=0.45\textwidth]{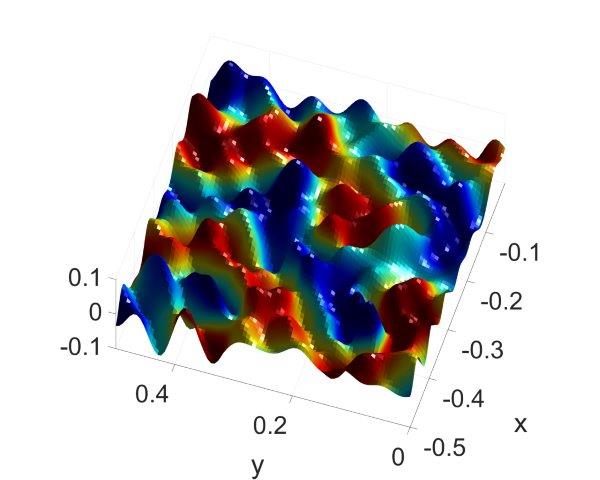}
		\end{overpic}
	\end{tabular}
	\caption{ Patterns  on rough surfaces {$\mathcal{M}$  with an} amplitude $\MAmp=0.1$ and $(M, N)\in \{(5, 15), (15,15)\}$. The model parameters for stripes are set according to Table~\ref{tbl:table1}.  Discretization parameters are $\ n_X=170=n_Y,  \tau=0.5,\ T=4000$.
}
	\label{fig:ex332}
	\end{figure}

These patterns, { both 2-D view and zoom in 3-D view}, are presented in {Figures}~\ref{fig:ex321}--\ref{fig:ex322} for spots, and in Figures~\ref{fig:ex331}--\ref{fig:ex332} for stripes. As before, the initial conditions are assigned to be random data in $[-0.5,0.5]$. In order to capture details of the spot and stripe patterns, we increase the spatial resolution from $n_X=90$ to $n_X=170$, again keeping $n_Y=n_X$ to ensure the irregular patterns {are} due to the roughness rather than low resolution.

Firstly, from Figures~\ref{fig:ex321} and \ref{fig:ex322}, we can see that the number of spots increases as the frequencies $M,N$ and the amplitude $\delta_M$ increase.
For {a} fixed amplitude $\MAmp$, patterns {undergo deformation} along the $x-$axis as the frequency $N$ increases. This {observation} is {particularly evident} in the case  $\MAmp=0.1$. On the other hand, for fixed values of $M,\ N$, the patterns maintain their {spot shapes}  when $\MAmp\leq 0.05$, and start to deform when $\MAmp \geq 0.05$. {When} $\delta=0.1$, all patterns become deformed spots.

Secondly, for {a} fixed $\delta_M=0.05$, zoom in 3-D profiles of region $[-0.3,0.1] \times [-0.3,0]$ for $[M,N]=[5,15]$ and region $[0.2,0.5] \times [-0.3,0]$ for $[M,N]=[15,15]$ reveal that irregular spots appear when part of the pattern is {situated} in a valley or ridge of the rough surface. When $\delta_M$ increases to $0.1$, as  in  Figure~\ref{fig:ex322}, the deformation of the patterns is more severe. {This is because the amplitude $\delta_M$ here is twice of that in  Figure~\ref{fig:ex321}}. From zoom in 3-D figures of region $[-0.1,0.5] \times [0.1,0.4]$ for $[M,N]=[5,15]$ and region $[-0.4,0.1] \times [-0.2,0.1]$ for $[M,N]=[15,15]$, the deformed patterns again appear when their locations cover the local ridge/valley/mountain of the rough surfaces.

{When} frequencies and amplitudes are varied, similar behavior (to the case of spots) can be observed in stripe formation, see {Figures}~\ref{fig:ex331} and \ref{fig:ex332}. High amplitudes and frequencies yield a particularly strong effect in the zoom in 3-D plots {displayed} in Figure \ref{fig:ex332}, where the stripes break into small separate components.  Further, we observe that the largest concentration values appear at the local peaks of the rough surface $\calm$.

\subsection{{Animal coat generation results}}\label{sec4:animal}
\begin{figure}[htbp]
\centering  
\label{animalcoatspatialfre}
\begin{tabular}{cccc}
 \begin{overpic}[width=0.45\textwidth]{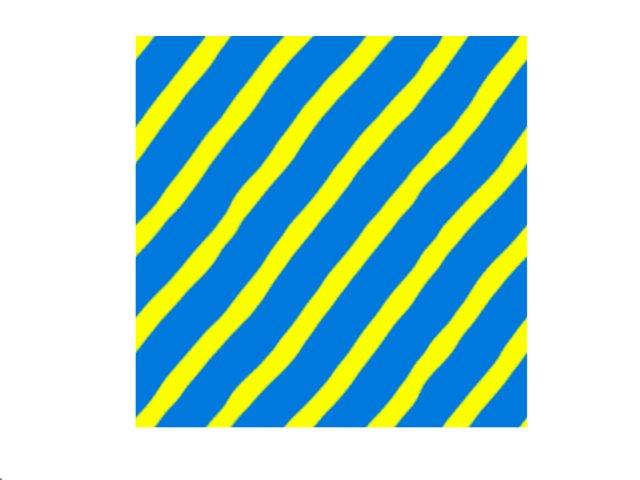}
 \put(20,0){(a) Emperor angelfish}
 \end{overpic}
 &
  \begin{overpic}[width=0.45\textwidth]{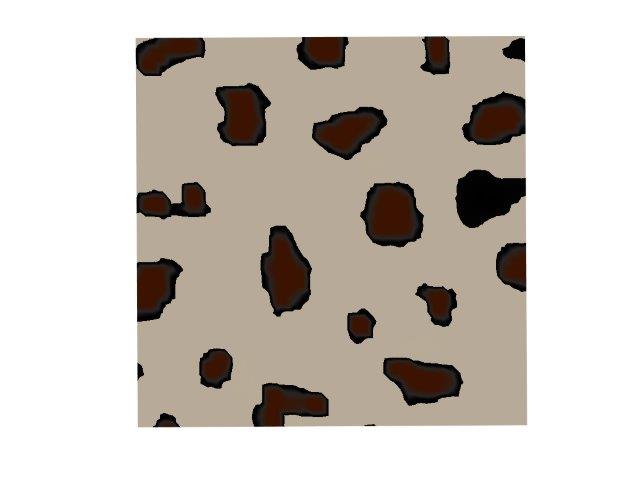}
 \put(40,0){(b) Genet}
 \end{overpic}
 \\
     \begin{overpic}[width=0.45\textwidth]{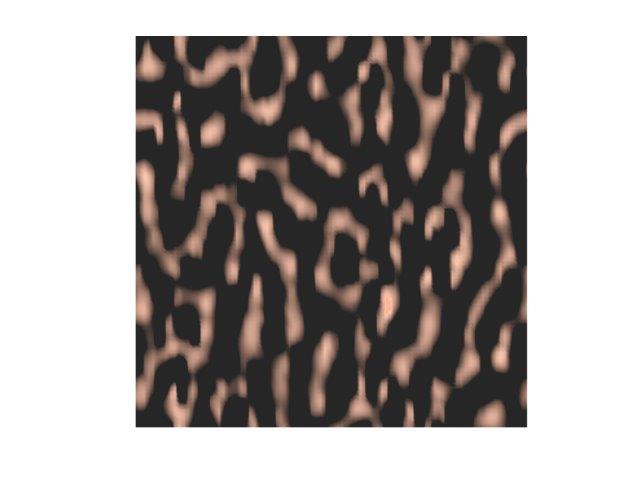}
     \put(30,0){(c) Plecostomi}
 \end{overpic}
 &
   \begin{overpic}[width=0.45\textwidth]{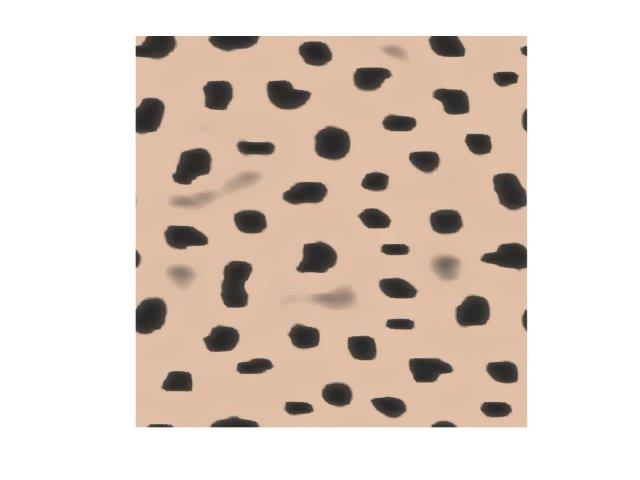}
 \put(40,0){(d) Cheetah}
 \end{overpic}
\end{tabular}
\caption{ {Experimental results for actual animal coat pattern formations on rough surfaces $\calm$ with parameters in Table~\ref{tbl:table1} and Table~\ref{tb:RSmanimal} }}\label{animal1}
\end{figure}

{When computing results on rough surfaces $\M$ characterized by spatial frequencies,} we observe interesting similarities between some steady state patterns and actual animal coat patterns. {In Figure~\ref{animal1},
we show the formation of the animal coats, as in Figure~\ref{fig:motivation}. These patterns are now on rough surfaces $\M$.}
The specific parameter values for generating each animal coat pattern are listed in Tables~\ref{tbl:table1} and \ref{tb:RSmanimal}.

\begin{table}[htbp]
	\centering
	\caption{\scriptsize  Parameters for {animal coats formation} on $\calm$ in Figure~\ref{animal1} with $n_X=n_Y=90$}
	\begin{tabular}{|l|c|c|c|c|c|c|}
	\hline
& & $M$  & $N$ &  $\tau$ & $T$ & $\delta_\mathcal{M}$  \\
\hline
Emperor angelfish \cite{Albert2014File}&Figure~\ref{animal1}~(a) & $5$  &$5$ & $0.5$ & $4000$ & $0.05$ \\
	\hline
Genet \cite{Category}& Figure~\ref{animal1}~(b) & $15$ &$15$ & $0.5$ & $800$ & $0.1$\\
\hline
Plecostomus \cite{Contributors2021Plecostomus} &Figure~\ref{animal1}~(c) & $15$ &$5$ & $0.5$ & $4000$ & $0.1$\\
\hline
Cheetah \cite{2015undefined} &Figure~\ref{animal1}~(d) & $15$ &$5$ & $0.5$ & $800$ & $0.1$\\
\hline
	\end{tabular}
	\label{tb:RSmanimal}
\end{table}

\section{Random rough surfaces $\cals$ by discrete data}\label{sec:RRS}
While the parametric equation \eref{eq:roughsurface1} {enables} us to work analytically on rough surfaces, i.e., via the evaluation of metric tensor \eref{eq:G}, its generalization to manifolds {\cite{Cheung+Ling-Kernembemethconv:18,Chen+Ling-Extrmeshcollmeth:20,ruuth_simple_2008}} is not  trivial.

In \cite{huang2021isogeometric}, the authors used the covariance function of random deformation fields and the surface
Karhunen-Lo\`{e}ve  expansion to generate random surfaces.  In \cite{hu1992simulation}, the authors {employed}  2D digital filters and Fourier analysis to generate random rough surfaces. In this section, we introduce a new approach for constructing rough surfaces $\cals$ based on random data and heat filters.
{Subsequently,} the desired reaction-diffusion systems are then solved on $\cals$.

\subsection{Construction of {$\cals$} by heat filters}

We {aim} to generate some random rough surfaces $\cals$, which have similar roughness as rough surfaces $\calm$ in \eref{eq:roughsurface} with different $M$ and $N$.

Let $[X,Y]\in \R^{n_X n_Y \times 2}\subset \cali^2$ as in Section~\ref{sec:FDM}. We assign uniform random numbers to each node to obtain the initial random surface values  $\tilde{Z}_0\sim \big(\calu[-1,1]\big)^{n_X n_Y }$ at nodes $[X,Y]$. We define the discretized heat filter according to the method {outlined in} Section~\ref{sec:FDM}, but with some filter-diffusion tensor  $\calf$ (instead of the diffusion tensor $\cala$ for surface $\calm$) in the Laplace-Beltrami operator $\triangle_{\calf,h} $. For an isotropic  filter, we {set} $\calf = I_{2\times 2}$.  This can generate an $\calm$-like surface with $M=N$. For the anisotropic case with $M\neq N$, we use $\calf = \text{diag}(2,1)$.
We can now \emph{smooth} the surface data $J$-times via
\[
    \tilde{Z}_{j+1} =(Q+\kappa h \triangle_{\calf,h})  \tilde{Z}_j, \qquad \text{for }j=0,\ldots, J,
\]
where $Q$ is an $n_X n_Y$ by $n_X n_Y$ matrix of ones, $\kappa>0$ is the parameter to control the weights of $\triangle_{\calf,h}$, and $h$ is the fill distance of the discrete set. This completes the definition of the pre-surface, that is the counterpart to $\tilde{Z}$ in
\eref{eq:roughsurface} of the surface type $\calm$. Similar to the scaling in \eref{eq:roughsurface}, we define $\cals$ {solely based on} nodal values
\[
    \Big \{ [X,Y, Z] :\, Z = \frac{\delta_\cals~}{\|\tilde{Z}_J\|_\infty} \tilde{Z}_J \Big\} \subset \cals,
\]
for some amplitude $\delta_\cals>0$.

\begin{figure}
	\centering
	\setlength{\abovecaptionskip}{0pt}
	\setlength{\belowcaptionskip}{0pt}
	\begin{tabular}{ccc}
	\begin{overpic}[width=0.3\textwidth]{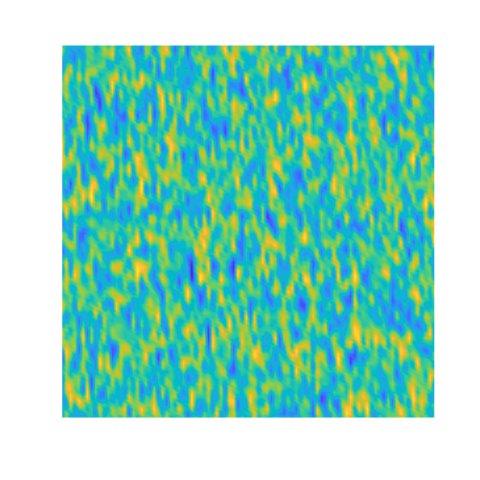}
		\put(20,5) {\scriptsize (a) $\calf=\text{diag}(0.01,1)$}
	\end{overpic}
	& \begin{overpic}[width=0.3\textwidth]{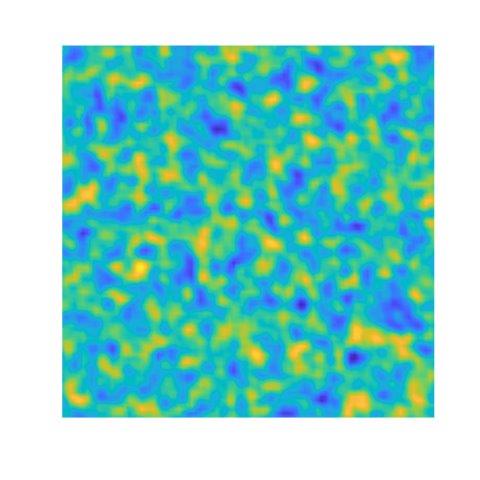}
		\put(20,5) {\scriptsize (b) $\calf=\text{diag}(1,1)$ }
		\end{overpic}
	& \begin{overpic}[width=0.3\textwidth]{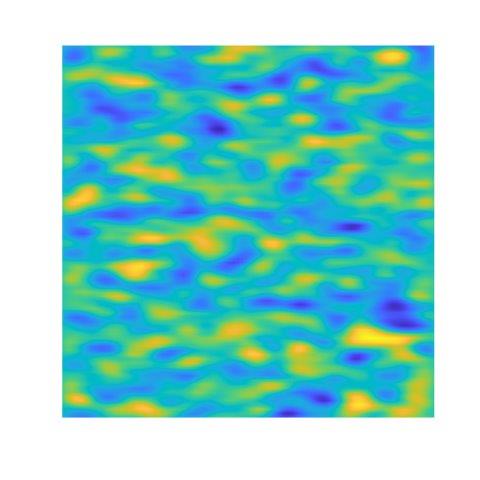}
	\put(25,5) {\scriptsize (c) $\calf=\text{diag}(5,1)$ }
\end{overpic}
	\end{tabular}
	\caption{ Rough surfaces $\mathcal{S}$ with $\kappa=2$ and the number of filter steps $J=15$,  {subject to different} filter-diffusion tensors $\calf$.
}
	\label{fig:rssdifffx}
	\end{figure}

{Subfigures of Figure~\ref{fig:rssdifffx} show} the rough surfaces $\cals$ derived for fixed $n_X=90,\ n_Y=n_X$, $\kappa=2$ and $J=15$ filter steps, {using different} filter-diffusion tensors
\[\calf \in \{\text{diag}(0.01,1), \text{diag}(1,1),\text{diag}(5,1)\},\]
respectively.
Taking $\calf(1,1)=1$ yields a rough surface $\cals$ which is similar to a rough surface $\calm$  with equal frequencies $M,N$ (cf. Figure~\ref{fig:roughsurface}).  {By setting $\calf(1,1)=0.01$, space is scales down in the $x$-direction,} yielding a rough surface $\cals$ which is similar to a rough surface $\calm$ with a larger frequency $M$, $M>N$. Conversely, {when} $\calf(1,1)=5$, a rough surface $\cals$ is obtained which {resembles} $\calm$ with a {lower} frequency $M, M<N$.

To reproduce the properties of surface $\calm$, the required value of  filter-diffusion tensors $\calf$ and number $J$ of filtering steps will depend on the density of the given data points $[X,Y]$.  For fixed amplitude $\delta_\mathcal{S}=1E-3$ and $n_X=90,\ n_Y=n_X$, Table \ref{tb:RSS} gives values of $\kappa,\  \calf$ and the number of filtering steps $J$ {required} to generate surfaces that {exhibit} good qualitative agreement with the rough surfaces $\calm$ in \eref{eq:roughsurface} for various $M,\ N$. Figure~\ref{fig:roughsurface} {provides} a comparison, {demonstrating} that rough surfaces $\cals$ by heat filters (column 2) are qualitatively similar to our {previous} surfaces  $\calm$ by \eref{eq:roughsurface} (see column 1).
\begin{table}[t]
	\centering
	\caption{\scriptsize  Coefficients for constructing the surfaces $\cals$  in the second column of Figure~\ref{fig:roughsurface}.  Parameters are chosen to give qualitative agreement with the rough surfaces $\calm$ \eref{eq:roughsurface} with $\delta_\mathcal{S}=1E-3,\ n_X=90=n_Y$}
	\begin{tabular}{|l|l|l|l|}
	\hline
	\multirow{3}{*}{$\calm$ in  \eref{eq:roughsurface} } & \multicolumn{3}{|c|}{$\cals$ by  heat filter } \\	
\cline{2-4}
& $\kappa$  & $\calf$ &  filter number     \\
\hline
 $[M,N]=[5,5]$& $5$ &$\text{diag}(1,1)$ & $15$   \\
	\hline
 $[M,N]=[5,15]$  &   $8$ & $\text{diag}(1,0.01)$  & $10$   \\
\hline
 $[M,N]=[15,15]$ & $0.2$ & $\text{diag}(20,20)$ & $2$    \\
\hline
	\end{tabular}
	\label{tb:RSS}
\end{table}

\begin{figure}
	\centering
	\setlength{\abovecaptionskip}{0pt}
	\setlength{\belowcaptionskip}{0pt}
	\begin{tabular}{ccc}
	\begin{overpic}[width=0.45\textwidth]{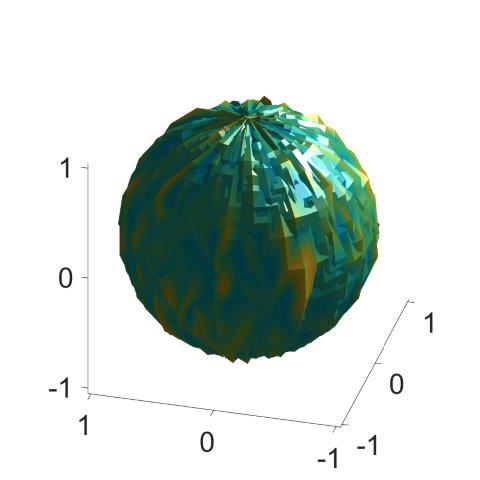}
		\put(40,0) { (a)}
	\end{overpic}
	& \begin{overpic}[width=0.45\textwidth]{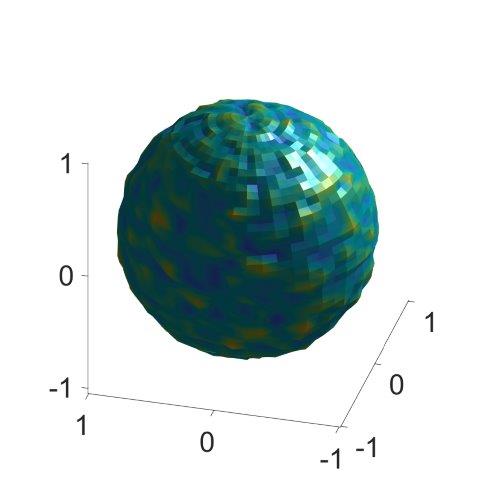}
		\put(40,0) { (b)}
		\end{overpic}
	\end{tabular}
	\caption{{Graphical illustration of rough closed manifolds by (a) type-$\calm$ procedure, where roughness is added to the (constant) radius function of the spherical coordinate, and  (b) type-$\cals$ procedure, where random noise is first added to data points on the unit sphere and then smoothed by heat filter.}}
	\label{fig:closedRoughSurface}
	\end{figure}

While not the focus of the current work, our construction methods for $\cals$ can be extended to generate rough closed manifolds {that are suitable for approximating numerical solutions of PDEs}. See Figure~\ref{fig:closedRoughSurface} for a graphical illustration of rough closed manifolds generated by type-$\calm$ and type-$\cals$ {construction methods}, respectively. Figure~\ref{fig:closedRoughSurface}(a) was obtained by applying the type-$\mathcal{M}$ procedure to the parameter space $(\theta,\phi)\in[0,2\pi]\times[0,\pi]$ with roughness added to the constant function $r=1$. Figure~\ref{fig:closedRoughSurface}(b) was obtained by adding noise to $r=1$ for points on the unit sphere. In contrast to  that in \eref{eq:disLB}, the heat filter here makes use of the discrete Laplace-Beltrami operator for the rough sphere.

\begin{figure}
	\centering
	\setlength{\abovecaptionskip}{0pt}
	\setlength{\belowcaptionskip}{0pt}
\begin{tabular}{cc}
\begin{overpic}[width=0.4\textwidth]{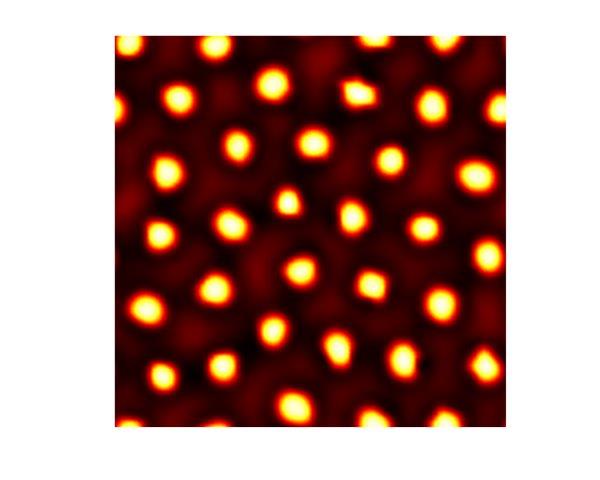}
	\end{overpic}
&
\begin{overpic}[width=0.4\textwidth]{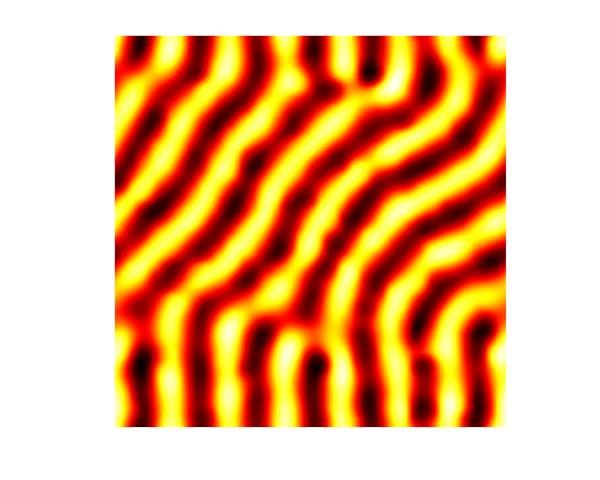}
	\end{overpic}
\end{tabular}
	\caption{For $\delta_\mathcal{S}=0.1,\ n_X=90, \ n_Y=n_X$, spot and stripe patterns on a rough surface $\cals$ with number of filter steps $J=15$, $\calf=\text{diag}(1,1)$, and $\kappa=5$.  These values from Table~\ref{tb:RSS} give a rough surface similar to $\calm$ with $M=5,N=15$.  Parameters for spots and stripes are {determined} according to Table~\ref{tbl:table1}.
}\label{fig:patternsRSS}
\end{figure}

The FDM in Section~\ref{sec:FDM} can be {employed} to solve PDEs on rough surfaces $\cals$. The only difference in the method is that we no longer have the parametric equation to calculate the metric tensor $G$ in \eref{eq:G}, and hence the diffusion tensor $\cala$ in \eref{eq:disLB}. Instead of computing metric tensor $G$ analytically as in {Sections}~\ref{sec:FDM} and \ref{sec:num},  
centered finite difference formulas are applied  to approximate the metric tensor  $G$. For solving reaction-diffusion systems on rough surfaces $\cals$, we set initial conditions to be steady state solutions on a zero amplitude rough surface $\calm$ \eref{eq:roughsurface}. Figure~\ref{fig:patternsRSS} plots the spot and stripe patterns on a rough surface $\cals$ using the reaction-diffusion parameters provided in Table~\ref{tbl:table1}. As {indicated} in the second line of Table~\ref{tb:RSS}, rough surface $\cals$ takes $\kappa=5$, $\calf=\text{diag}(1,1)$ with $J=15$ filter steps  to approximate rough surface $\calm$ with $M=5,N=5$.  It {is evident} that the spot and stripe patterns generated on $\cals$ are {closely resemble} those on $\calm$ under parameters $M=5,N=5,\ \delta_\mathcal{M}=0.1$ (see the third row of Figure~\ref{fig:ex31}).

\begin{figure}[htbp]
\centering  
\begin{tabular}{cccc}
\begin{overpic}[width=0.45\textwidth]{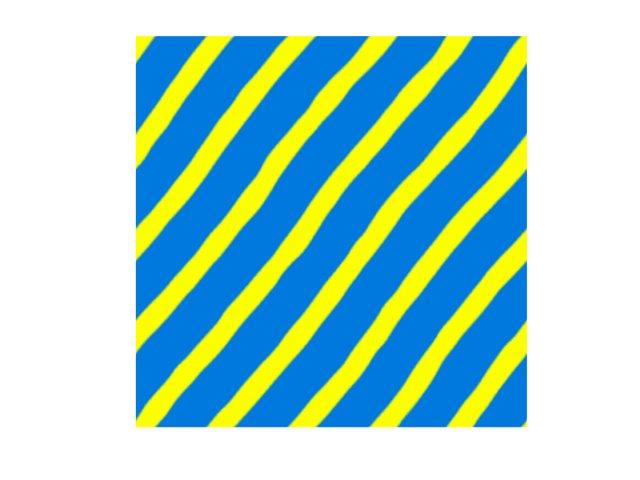}
\put(20,0){(a) Emperor angelfish}
\end{overpic}
&
\begin{overpic}[width=0.41\textwidth]{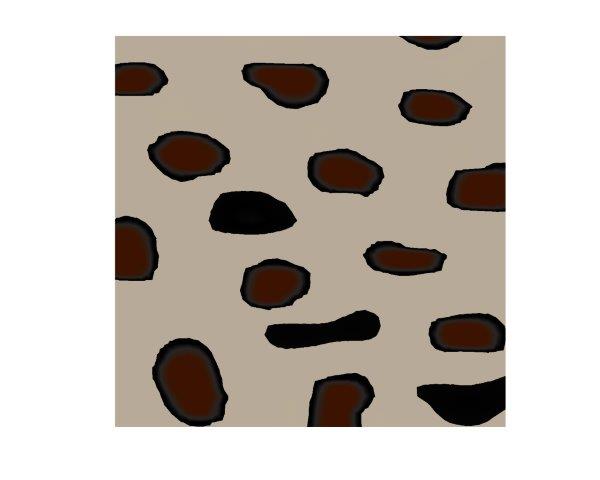}
\put(35,0){(b) Genet }
\end{overpic}
\\
\begin{overpic}[width=0.45\textwidth]{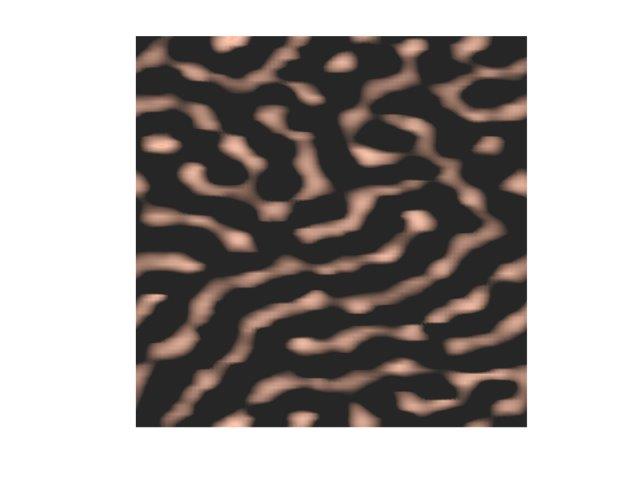}
\put(30,0){(c) Plecostomi}
\end{overpic}
&
\begin{overpic}[width=0.42\textwidth]{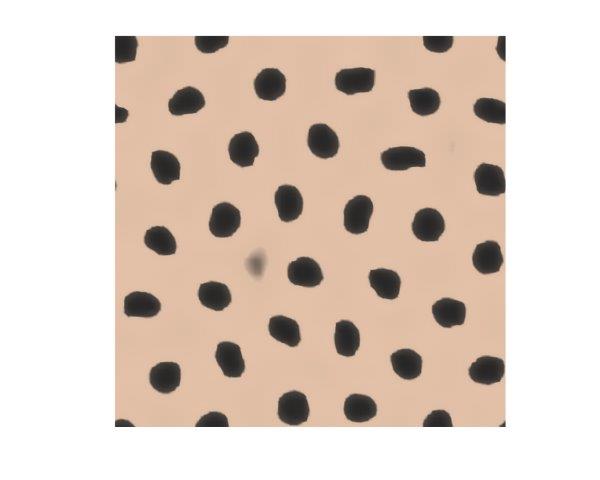}
\put(35,0){(d) Cheetah}
\end{overpic}
\end{tabular}
\caption{{Experimental results for actual animal coat pattern formations on rough surfaces $\cals$ with parameters in Table~\ref{tbl:table1} and Table~\ref{tb:RSSmanimal}. }}
\label{animalcoatSsurface}
\end{figure}

We conclude by generating the same set of animal coats displayed in {Figure~\ref{fig:motivation}} using Turing models on rough surfaces $\mathcal{S}$ with parameters {specified in} Table~\ref{tbl:table1} and  Table~\ref{tb:RSSmanimal}. Figure~\ref{animalcoatSsurface} demonstrates again that adding surface roughness into the process of pattern generation can indeed provide more varied formations for animal coats. Moreover, the patterns on surface $\cals$ {exhibit} valid formations just like $\M$.

\begin{table}[htbp]
	\centering
	\caption{\scriptsize  Parameters for {animal coat generation} on $\mathcal{S}$ with $n_X=n_Y=90$}
	\resizebox{\textwidth}{12mm}{
	\begin{tabular}{|l|c|c|c|c|c|c|c|}
	\hline
& & $\kappa$  & $\mathcal{F}$ & $J$ &  $\tau$ & $T$ & $\delta_\mathcal{S}$  \\
\hline
Emperor angelfish  \cite{Albert2014File}&Figure~\ref{animalcoatSsurface}~(a) & $5$ & $\text{diag}(1,1)$ & $15$ & $0.5$ & $400$ & $0.05$\\
	\hline
Genet\cite{Category} &Figure~\ref{animalcoatSsurface}~(b) & $8$ &$\text{diag}(1,0.01)$ & $10$ & $0.5$ & $800$ & $0.1$\\
\hline
Plecostomus \cite{Contributors2021Plecostomus} &Figure~\ref{animalcoatSsurface}~(c) & $8$ &$\text{diag}(1,0.01)$ & $10$ & $0.5$ & $3000$& $0.1$\\
\hline
Cheetah \cite{2015undefined}& Figure~\ref{animalcoatSsurface}~(d) & $0.2$ &$\text{diag}(20,20)$ & $2$ & $0.5$ & $400$& $0.05$\\
\hline
	\end{tabular}}
	\label{tb:RSSmanimal}
\end{table}

\section{Conclusion} In this paper, we discussed the generation of surfaces with arbitrary roughness and the formation of patterns on such surfaces. Rough surfaces are characterized by two different methods: surfaces by analytic parametric equations, and surfaces constructed by imposing discretized heat filters on random nodal values.  The patterns generated on random rough surfaces {with varying} amplitudes and spatial frequencies are illustrated, and the computational results indicate that the patterns became {increasingly} irregular as amplitude increases. The change of spatial frequencies on the $x$ and $y$ axes  also lead to pattern deformation along the {corresponding} directions.  It can be concluded that surface roughness leads to the patterns becoming deformed and of different sizes. Conversely, when the amplitude equals zero, i.e., a flat two dimensional domain, the patterns are regular stripes and spots.  We conclude that combining reaction-diffusion systems
with rough surfaces gives a way to {achieve} a better match to
the observed variability of patterns in a variety of real world animal coats. Furthermore, the method for generating rough surfaces by  heat filters can be further applied to obtain closed rough manifolds. We plan to explore this generalization in future work.


\section*{Appendix: Finite difference algorithm for the heat equation}
The heat equation \eref{eq heat}
on a rough surface defined over the parameter space $\mathcal{I}^2=[-1,1]^2$ is considered.
The set of discrete data points on $\mathcal{I}^2$ {is} defined as:
$$[X,Y]:=\big\{\big\{(x^i_1,x^j_2)\big\}_{i=1}^{n_X}\big \}_{j=1}^{n_Y} \in \R^{n_X n_Y \times 2 },$$
with mesh size $h_{x_1},h_{x_2}$ on each axis.
To fully discretize \eref{eq:timeBE}, we require
discretization of the Laplacian-Beltrami operator. For
any twice differentiable function $u:\cali^2\to\R$, we introduce differentiation matrices $\mathcal{D}_k \in \R^{n_X n_Y \times n_X n_Y },\ k\in \{1,2\} $  satisfying periodic boundary condition \eref{eq:HeatBC}
as follows. Using a second-order centered finite differences, we have
for each point $(x^i_1,x^j_2)$ that
\begin{equation*}
\begin{aligned}
 \frac{\partial u}{\partial x_1}(x^i_1,x^j_2)
 =\begin{bmatrix}
 -\frac{1}{2h_{x_1}} &\frac{1}{2h_{x_1}}
\end{bmatrix}
\begin{bmatrix} u(x^{i-1}_1,x^j_2) \\ u(x^{i+1}_1,x^j_2)\end{bmatrix} , \\
 \frac{\partial u}{\partial x_2}(x^i_1,x^j_2)
 =\begin{bmatrix}
 -\frac{1}{2h_{x_2}} &\frac{1}{2h_{x_2}}
\end{bmatrix}
\begin{bmatrix} u(x^{i}_1,x^{j-1}_2) \\ u(x^{i}_1,x^{j+1}_2)\end{bmatrix}.
\end{aligned}
\end{equation*}
The following fictitious node approach is applied to {handle} the periodic boundary conditions:
\begin{equation*}
\begin{aligned}
 &u(x_1^0,x_2^j)=u(x_1^{n_X-1},x_2^j), \ u(x_1^{n_X+1},x_2^j)=u(x_1^2,x_2^j),\ u(x_1^i,x_2^0)=u(x_1^i,x_2^{n_Y-1}),\\
 &u(x_1^i,x_2^{n_Y+1})=u(x_1^i,x_2^2),\
 u(x_1^1,x_2^j)=u(x_1^{n_X},x_2^j), \quad \ u(x_1^i,x_2^1)=u(x_1^i,x_2^{n_Y}).
\end{aligned}
\end{equation*}
By assembling all points in $[X,Y]$,
the nodal values of $\frac{\partial u}{\partial x_k}$ ($k=1,2$) at $[X,Y]$ can be obtained by
\[
  \frac{\partial u}{\partial x_k}(X,Y) \approx \mathcal{D}_k u(X,Y), \qquad \text{for } k\in\{1,2\},
\]
where
$u(X,Y):= [u(x_1^1,x_2^1),\cdots,u(x_1^{n_X},x_2^1), \cdots,u(x_1^1,x_2^{n_Y}),\cdots,u(x_1^{n_X},x_2^{n_Y})]^{T}$
is the vector of nodal function values and differential matrices $\mathcal{D}_k \in \mathbb{R}^{n_Xn_Y \times n_Xn_Y}$ are in the form of
\begin{equation*}
\begin{aligned}
\mathcal{D}_{1}&=
 \begin{bmatrix}
  D_{1}^{B} & \cdots&0    \\
  \vdots&  \ddots&   \vdots\\
     0&   \cdots &D_{1}^{B}
\end{bmatrix}, \
D_{1}^{B}=\frac{1}{2h_{x_1}}
 \begin{bmatrix}
  0& 1 & 0&\cdots& 0 &-1 & 0   \\
  -1& 0 & 1 &\cdots& 0 & 0 & 0  \\
  \vdots& \vdots&\vdots & \ddots&  \vdots& \vdots  \\
     0& 0& 0 &\cdots& -1 & 0 & 1\\
  1& 0 & 0&\cdots& 0 &0 & -1
\end{bmatrix},
\end{aligned}
\end{equation*}
and
\begin{equation*}
\begin{aligned}
\mathcal{D}_{2}=
 \begin{bmatrix}
   0 &D_{2}^{B} & 0& \cdots &0 &-D_{2}^{B} &0    \\
  -D_{2}^{B} & 0 & D_{2}^{B}& \cdots& 0 & 0 &0 \\
  \vdots& \vdots& \vdots& \ddots& \vdots &  \vdots &   \vdots\\
0& 0&  0& \cdots &-D_{2}^{B}&0& D_{2}^{B}\\
     D_{2}^{B} & 0&  0& \cdots &0 &   0&- D_{2}^{B}
\end{bmatrix}, \ D_{2}^{B}=\frac{1}{2h_{x_2}}
 \begin{bmatrix}
   1 &\cdots&0 \\
  \vdots&\ddots& \vdots \\
   0 &\cdots&1
\end{bmatrix},
\end{aligned}
\end{equation*}
with $D_{k}^{B} \in \R^{n_X\times n_X},\ k\in \{1,2\}$.

\bibliographystyle{siam}
\bibliography{RS_FDM}


\end{document}